\documentclass{amsart}
\usepackage{amsfonts}
\usepackage{amsmath}
\usepackage{amssymb}
\usepackage{graphicx,array}

\newtheorem{theorem}{Theorem}[section]
\theoremstyle{plain}
\newtheorem{definition}[theorem]{Definition}
\newtheorem{corollary}[theorem]{Corollary}
\newtheorem{lemma}[theorem]{Lemma}

\def\sideremark#1{\ifvmode\leavevmode\fi\vadjust{\vbox to0pt{\vss
 \hbox to 0pt{\hskip\hsize\hskip1em
 \vbox{\hsize2cm\tiny\raggedright\pretolerance10000
 \noindent #1\hfill}\hss}\vbox to8pt{\vfil}\vss}}}%

\newcommand\straightparen[1]{{\rm(}{#1}{\rm\,)}}
\newcommand\usfrac[2]{{#1}/{#2}}
\renewcommand{\arraystretch}{1.3}
\def\dint{{\displaystyle\int}}

\numberwithin{equation}{section}

\numberwithin{equation}{section}

\newcommand{\cF}{\mathcal{F}}

\newcommand{\cU}{\mathcal{U}}
\newcommand{\cV}{\mathcal{V}}

\newcommand{\GL}{\mathrm{GL}}

\newcommand{\fa}{\mathfrak{a}}
\newcommand{\fb}{\mathfrak{b}}
\newcommand{\fc}{\mathfrak{c}}
\newcommand{\fd}{\mathfrak{d}}
\newcommand{\fe}{\mathfrak{e}}
\newcommand{\ff}{\mathfrak{f}}
\newcommand{\fg}{\mathfrak{g}}

\newcommand{\fk}{\mathfrak{k}}

\newcommand{\fn}{\mathfrak{n}}

\newcommand{\fp}{\mathfrak{p}}
\newcommand{\fq}{\mathfrak{q}}

\newcommand{\fu}{\mathfrak{u}}
\newcommand{\fz}{\mathfrak{z}}

\newcommand{\ga}{\alpha}

\newcommand{\gD}{\Delta}
\newcommand{\gG}{\Gamma}
\newcommand{\gL}{\Lambda}

\newcommand{\R}{\mathbb{R}}
\newcommand{\C}{\mathbb{C}}
\newcommand{\N}{\mathbb{N}}
\newcommand{\Z}{\mathbb{Z}}
\newcommand{\bc}{\mathbf{c}}

\newcommand{\Supp}{\mathrm{Supp}}
\newcommand{\Ad}{\mathrm{Ad}}
\newcommand{\re}{\mathop{\mathrm{Re}}}
\newcommand{\im}{\mathop{\mathrm{Im}}}
\newcommand{\id}{\mathrm{id}}
\newcommand{\Exp}{\mathrm{Exp}}
\newcommand{\rank}{\mathrm{rank}}

\newcommand{\e}{\varepsilon}
\renewcommand{\l}{\lambda}

\newcommand{\inner}[2]{\langle #1,#2\rangle}

\begin{document}
\makeatletter
\title[Paley-Wiener Theorem]{The Paley-Wiener Theorem and the local Huygens'
principle\\ for compact symmetric spaces}
\author{Thomas Branson}
\address{Department of Mathematics, University of Iowa, 
Iowa City IA 52242 USA}
\email{thomas-branson@uiowa.edu}
\urladdr{http://www.math.uiowa.edu/\symbol{126}branson}
\author{Gestur \'{O}lafsson}
\address{Department of Mathematics, Louisiana State University, 
Baton Rouge LA 70803, USA}
\email{olafsson@math.lsu.edu}
\urladdr{http://www.math.lsu.edu/\symbol{126}olafsson}
\author{Angela Pasquale}
\address{Laboratoire et D\'epartement de Math\'ematiques,
  Universit\'e de Metz, France}
\email{pasquale@math.univ-metz.fr}
\urladdr{http://www.math.univ-metz.fr/\symbol{126}pasquale}
\thanks{GO was supported by NSF grants DMS-0070607, DMS-0139783,
and DMS-0402068,  
and by the DFG-Schwerpunkt ``Global Methods in Complex Geometry.''
TB was partially supported by the Erwin Schr\"odinger Institute.}
\date{November 13, 2004}
\subjclass[2000]{Primary 53C35; Secondary 35L05}
\keywords{}
\dedicatory{Dedicated to Gerrit van Dijk on the occasion of his 65th birthday}

\begin{abstract} We prove a Paley-Wiener Theorem for 
a class of symmetric spaces of the 
compact type, in which all root multiplicities are even.
This theorem characterizes functions of small support in terms of holomorphic
extendability and exponential type of their (discrete) Fourier transforms.
We also provide three independent new proofs of the strong Huygens' principle
for a suitable constant shift of the wave equation on odd-dimensional spaces
from our class.

\end{abstract}
\maketitle
\makeatother

\section*{Introduction}

\noindent
In the context of spherical harmonic analysis, the compactness of
a symmetric space $U/K$ is reflected by the discreteness of its dual space,
which is the set of irreducible $K$-spherical unitary
representations of $U$. The same set
parameterizes the set of (elementary) spherical functions.
Thus the spherical Fourier transforms of $K$-invariant
functions on $U/K$ are functions on a discrete set.
Likewise, the formula for spherical inversion, which recovers a sufficiently
regular function on the symmetric space in terms spherical functions, 
is given by a series.
This structural discreteness can be overcome for functions with 
``small support'', by relating them to functions on the
tangent space $\mathfrak{q}$ to $U/K$ at the base point $x_0=\{K\}$.
This procedure can be easily illustrated in the Euclidean setting:
consider a smooth function $f:S^1 \to \C$, where $S^1$ denotes the
unit circle. View $f$ as a periodic functions on $\R$ by
$t\mapsto f(e^{it})$ and assume  that $f$ has small support, say
in  $[-R,R]+2\pi \Z$, where $0<R<\pi$.
We can then regard $f$ as a smooth function on the real line with support
in $[-R,R]$ by setting it equal to $0$ outside of the fundamental
period $[-\pi,\pi)$.
By the classical Paley-Wiener theorem on $\R$, the Fourier transform of
$f$ is an entire function of exponential type $R$. 
It therefore provides a holomorphic extension of the Fourier transform of
$f$ as a function on $S^1$.
Likewise, if $F$ is a holomorphic function
on $\C$ of exponential type $R$, $0<R<\pi$,
then the inversion formula for the continuous Fourier transform
gives a function $f_1$ with support in $[-R,R]$, and we
can define a function $f$ on $S^1$ by $f(e^{it})=f_1(t)$.

The possibility of characterizing central smooth functions with ``small 
support'' on compact Lie groups by means of the entire extension 
and exponential growth of their Fourier transform was first proven by 
Gonzalez in \cite{Go00}.
In this paper we extend the local Paley-Wiener theorem to all 
compact symmetric spaces $U/K$ with even multiplicities:
the $K$-invariant smooth functions on $U/K$ with ``small support'' 
will be characterized in terms of holomorphic extendibility
and exponential growth of their spherical Fourier transform.
Moreover, the exponential growth of the transformed function will
be linked to the size of the support of the function on the symmetric space.
The given characterization relies on the fact that the spherical
functions on a compact symmetric space extend holomorphically to
the complexified symmetric space. Their restrictions to the
noncompact dual symmetric spaces $G/K$ are in turn spherical functions
on $G/K$. This allows us to use
known information on the
spherical functions on $G/K$
and classical Fourier analysis on the
Lie algebra of a maximal abelian
subspace of $\fq$.
In particular, the classical Paley-Wiener theorem is used to obtain the
required holomorphic extension of compact spherical Fourier
transforms of a $K$-invariant function on $U/K$ with a small support.

Properties of holomorphic extendibility for spherical functions
on symmetric spaces have been the 
objects of intensive recent study, with different approaches and perspectives.
See e.g.\ \cite{KS1}, \cite{Pthesis}, and \cite{Opd94}.
The situation which we consider in this paper corresponds to
symmetric spaces with even multiplicities.
It is rather special because of the existence of shift operators providing
explicit formulas for the spherical functions by relating
them to exponential functions \cite{OP04}.
These shift operators are suitable multiples of Opdam's shift operators.
The multiplying factor has been chosen so, as to cancel the singularities of
the
coefficients of Opdam's shift operators along the walls of the Weyl
chambers. The resulting operators are differential operators with
holomorphic coefficients. Hence, we can read off
the properties of holomorphic extendibility of the spherical functions
directly from these formulas.
Furthermore, the shift operators allow us, as mentioned above, to reduce
several problems in harmonic analysis on symmetric spaces of even 
multiplicities
to the corresponding problems in Euclidean
harmonic analysis.

Our proof depends heavily on the assumption that
all root multiplicities are even, and it is not possible to 
generalize it to obtain local Paley-Wiener type theorems
for general compact symmetric spaces. On the other hand, the same proof 
can be employed for several other even-valued multiplicity
functions which are not geometric. We will not explore this avenue 
further in the present article.

The relation between spherical transforms on compact and noncompact
symmetric spaces investigated in this paper also yields a representation of
smooth functions with ``small support'' on the compact space as
integrals of spherical functions of the noncompact dual.
These integral formulas are the key ingredient for studying the solutions
of the wave equation on Riemannian symmetric spaces of compact type.
From exponential estimates for the solutions, we deduce in Section
\ref{section:Huygens} that the strong
Huygens' principle is valid on these spaces.

The \textit{\straightparen{strong} Huygens' principle} states that,
in odd dimensions,
the light at time $t_0$ at a location $x$ influences at later times $t_1$ only
those locations which have distance exactly $t_1-t_0$ from $x$.
Hence,
if a wave is supported in the sphere $\{x\mid\|x\| \leq  R \}$
at the initial time $0$, then it will be supported in the annulus
$\{x\mid t-R \leq \|x\|\leq t+R\}$ at time $t$.
In particular, at times $t >R$, the wave will vanish inside the
sphere $\{x\mid \|x\| < t-R \}$.

Several different authors have proven the validity of Huygens' principle on
odd dimensional Riemannian symmetric spaces with even multiplicities of
either noncompact or compact type.  Here ``light'' is to be interpreted
as a solution of a suitable wave equation, obtained by a certain constant
shift
of the d'Alembertian.
Their proofs use a variety of different methods.
The first results in this direction were given by Helgason \cite{He64,He84a},
see also \cite{He94},
who proved Huygens' principle for symmetric spaces $G/K$ for which either
$G$ is complex or $G={\rm SO}_0(n,1)$, and for compact groups.
In the general case of odd dimensional
Riemannian symmetric spaces of the noncompact type with
even multiplicities,
the validity of Huygens' principle was stated
without proof by Solomantina \cite{So}.
A proof by Radon transform methods was provided by \'Olafsson and 
Schlichtkrull \cite{OS92}.
An independent proof was obtained by Helgason \cite{He92} by means
of his Fourier transform.
In \cite{BOS95} the authors proved an exponential decay property for
solutions of the wave equation with compactly supported
initial data. This method implied another
independent proof of the Huygens' principle
for odd dimensional symmetric spaces with
even multiplicities; see \cite{BOS95}.
Finally, a completely different approach based on Heckman-Opdam's shift
operators and explicit formulas for the fundamental solutions
was provided by Chalykh and Veselov in \cite{CV96}.
The formulas of Chalykh and Veselov give the fundamental solution of the 
wave equation in polar coordinates.
By replacing hyperbolic functions with their trigonometric counterparts,
one can also deduce formulas for the fundamental solutions of the
wave equation on compact symmetric spaces. These
formulas will be valid for small values of time. Using this argument,
Chalykh and Veselov state that Huygens' principle holds also on
Riemannian symmetric spaces of the compact type.

In the context of Riemannian symmetric spaces,
Huygens' principle has been much less studied for compact type
than for noncompact type.
In \cite{bent}, {\O}rsted used conformal properties of wave operators and
of Lorentzian spaces covered by $\R\times S^{2n+1}$ to establish Huygens'
principle for the wave, Dirac, and Maxwell equations on $S^{2n+1}$.  
His proof makes it
clear that analogues will be valid for other linear differential operators
with suitable hyperbolicity and conformal properties.
A different 
proof for the wave equation on the odd sphere 
$S^{2n+1}$ were given by Lax and Phillips \cite{LP78}.
Branson \cite{tbjfa} extended the Lax-Phillips proof to an infinite class
of hyperbolic equations on the odd sphere.
Helgason proved Huygens' principle for the compact group case, see
\cite{He84a}.
Finally, Branson and \'Olafsson \cite{BO97} proved
that the local Huygens' principle for a compact symmetric space
$U/K$ is valid if and only if
Huygens' principle holds for the non-compact dual space
$G/K$.

In this article we provide three independent new proofs of a local version
of the strong Huygens' principle for compact symmetric spaces $U/K$. 
One of these methods comes from exponential 
estimates for the smooth solutions of $K$-invariant Cauchy problems 
for the modified wave equations on $U/K$. These estimates are 
obtained by methods similar to those introduced for the noncompact 
setting in \cite{BOS95}. 
It is nevertheless important to mention that the use of the shift operators
indeed reduces the proof of the of Huygens' principle 
on Riemannian symmetric
spaces of either type (compact or noncompact) to
the validity of the same principle in the Euclidean setting. The proof
presented in this paper is therefore easier than that in \cite{BOS95}.

Another proof of the local strong
Huygens' principle is in the spirit of the paper of Chalykh and
Veselov \cite{CV96}. The formulas for the spherical functions proven
in Theorem \ref{th-extension} permit us to derive an explicit formula for
the solution
of the wave equation corresponding to a given smooth initial condition.
The tools for writing down these formulas appear in the proof of the
local Paley-Wiener theorem.  An essential property in our argument is 
that our shift operators, which link spherical functions to
exponential functions, have regular (indeed analytic) coefficients. 
This fact was not proven in \cite{CV96}.
  
Our paper is organized as follows. In Section \ref{section:symmspaces} 
we recall some structure theory of Riemannian symmetric spaces of the
compact type. The spherical functions and spherical representations 
are introduced in Section \ref{section:spherical}. Theorem \ref{th-extension}
proves the existence of differential shift
operators. These provide explicit formulas for the spherical functions 
on compact symmetric spaces. The main theorem in this paper 
is the local Paley-Wiener theorem, which is stated and proven in 
Section \ref{section:localPW}. The integral formula for functions
with ``small support'' is given by Corollary \ref{cor:integralformula}.
Finally, Section \ref{section:Huygens} contains the proofs of the local 
strong Huygens' principle on Riemannian symmetric spaces of the compact type.

\section{Symmetric spaces}
\label{section:symmspaces}

\subsection{Compact symmetric spaces}

\noindent
In this section we recall some facts about compact symmetric
spaces. We use \cite{He78}, Chapter VII, and \cite{Take}, Chapter II, 
as standard references.

Let $U$ be a connected compact Lie group with center $Z$ and
Lie algebra $\fu$. Denote by $\fz$ the center of $\fu$. Then
$\fu =\fz\oplus \fu^\prime$, where $\fu^\prime:=[\fu,\fu ]$ 
is semisimple.
Let $\exp:\fu \to U$ be the exponential map.
If 
$\fz\not=\{0\}$, then we set
$\Gamma_0:=\{X \in \fz\mid \exp X=e\}$, where $e$ denotes the identity of
$U$. Then $\Gamma_0$ is a full rank lattice in $\fz$ and
$T:=\fz/\Gamma_0$ is isomorphic to the identity connected component of $Z$.
We will from now on write $T=Z_0$.
Denote by $U^\prime$ the analytic subgroup of $U$ with
Lie algebra $\fu^\prime$. 
Then $U^\prime$ is semisimple with finite center and
$U=TU^\prime\simeq T\times_F U^\prime$
where $F=T\cap U^\prime$ is a finite
central subgroup of $U$.
We will for simplicity assume that $F$ is trivial.
Thus $U\simeq T\times U^\prime$.

Let $\tau:U\rightarrow U$ be a non-trivial analytic
involution.
Set $U^\tau :=\{u\in U\mid \tau (u)=u\}$, and define
$K$ be the identity connected component of $U^\tau$.
Then $U/K$ is a connected  compact symmetric space (also called
Riemannian symmetric space of the compact type).
The derived involution of $\tau$ on $\fu$ will be denoted by the same 
letter $\tau$. Thus $\tau(\exp(X))=\exp(\tau(X))$ for all $X\in\fu$.

Let $\fk$ denote the Lie algebra of $K$.
We shall assume that $U$ acts effectively on $U/K$, i.e. that 
$\fk \cap \fz=\{0\}$. 
 Then
\begin{equation*}
\mathfrak{k}=\fu^\tau :=\{X\in\mathfrak{u}\mid\tau(X)=X\}\subset\fu^\prime\, .
\end{equation*}
Set
$$\mathfrak{q}:=\{X\in\mathfrak{u}\mid\tau(X)=-X\}\, .$$
Then $\mathfrak{u}=\mathfrak{k}\oplus\mathfrak{q}$ and
$\fz\subseteq \fq$.

For a real vector space $V$ we denote by $V^\ast$ its dual and by
 $V_{\mathbb{C}}:=V\otimes_{\mathbb{R}}\mathbb{C}$ its complexification.
If $V$ is a Euclidean vector space with inner product
$\inner{\cdot}{\cdot}$ and $W\subseteq V$ is a subspace,
then $W^{\perp}$ denotes the orthogonal complement of
$W$ in $V$. We identify
$W^\ast$ with the space $\{f\in V^{\ast}\mid f|_{W^{\perp}}=0\}$.
The complex linear extension to $V_\C$ of a linear map $\varphi:V\to V$ 
will be denoted by the same symbol $\varphi$.
For $\lambda\in V^*$  define $h_\lambda \in V$ by
$\lambda (H)=\inner{H}{h_\lambda}$. For $\lambda \not= 0$ we
set $H_\lambda :=2\inner{h_\lambda}{h_\lambda}^{-1}h_\lambda$. Then
$\lambda (H_\lambda)=2$. Finally we define an inner product on $V^\ast$
by $$\inner{\lambda}{\mu}:=
\inner{h_\lambda}{h_\mu}=\lambda (h_\mu)=\mu (h_\lambda)\,.$$

Recall that the Killing form $\kappa$ on $\fu$ is negative definite 
on $\mathfrak{u}^{\prime}$. 
Fix an inner product $\inner{\cdot}{\cdot}$ on $\fz$ and define a 
$U$-invariant inner product on $\mathfrak{u}$ by
\begin{equation*}
\inner{Z_{1}+X_{1}}{Z_{2}+X_{2}}:=\inner{Z_{1}}{Z_{2}}
-\kappa(X_{1},X_{2})\,, \qquad
Z_{1},Z_{2}\in \fz\, ,\; X_{1},X_{2}\in\mathfrak{u}^{\prime}\, .
\end{equation*}

Let $\mathfrak{b}\subseteq\mathfrak{q}$ be a
maximal abelian subspace and set $\fb_1:=\fb\cap \fu^\prime$. Then
$$\mathfrak{b}=\fz \oplus \mathfrak{b}_1\, .$$
Set $\mathfrak{a}:=i\mathfrak{b} \subseteq \mathfrak{u}_{\mathbb{C}}$
and $\fa_1=i\fb_1$. 
Then, by restriction, $\inner{\cdot}{\cdot}$ defines an inner product on
$\fa$, and hence we can apply the above notational conventions to $(\fa,
\inner{\cdot}{\cdot})$.
In particular $H_\lambda\in \fa$ is well defined for all nonzero 
$\lambda\in \fa^*$.

For $\alpha\in\mathfrak{b}_{\mathbb{C}}^{\ast}=\fa^*_\C$ let
\begin{equation*}
\mathfrak{u}_{\mathbb{C}}^{\alpha}:=\left\{
X\in\mathfrak{u}_{\mathbb{C}}\mid\forall H\in\mathfrak{b}\, :\,  [H,X]=\alpha(H)X\right\}
\end{equation*}
and set $m_{\alpha}:=\dim_{\mathbb{C}}\mathfrak{u}_{\mathbb{C}}^{\alpha}$.
If $\mathfrak{u}_{\mathbb{C}}^{\alpha}\not =\{0\}$, then 
$\ga$ is called a \textit{root} and $m_\alpha$ is its \textit{multiplicity}.
We denote by $\Delta$ the set of roots and by $W=W(\Delta)$ the
corresponding Weyl group. Recall that $W$ is generated by the
reflections $s_\alpha$ with $\alpha \in\Delta$. 
Here $s_\alpha(H):= H-\alpha(H)H_\alpha$.
If $\ga\in \Delta$, then $\mathfrak{u}_{\mathbb{C}}^{\alpha}
\subseteq\mathfrak{u}_{\mathbb{C}}^{\prime}$,
$\alpha |_{\fz_\C}=0$, 
and $\alpha\in i\mathfrak{b}^{\ast}_1=\mathfrak{a}^{\ast}_1$. 
Hence $\ga$ is real valued on $\fa$ and $\ga|_{\fz}=0$.
Choose $X\in \fa$ so that $\alpha (X)\not= 0$ for
all roots $\alpha$.  Then $\Delta^{+}:=\{\alpha \in \Delta\mid
\alpha (X)>0\}$ is a set of \textit{positive roots}. We denote by 
$\Sigma$ the corresponding set of simple roots.

\subsection{Integration on $U/K$}

We now fix our normalization of measures. If $L$ is
a locally compact Hausdorff topological group, then
$dl$ denotes a left invariant (Haar) measure
on $L$.
When $L$ is a compact group we normalize $dl$ 
so that the volume of $L$ is $1$. In this case, if $M$ is a closed 
(and hence compact) subgroup of $L$, 
then we normalize the invariant measure $d(lM)$ on $L/M$
so that $L/M$ has volume $1$. 
We then have for all $f\in  L^1(L/M)$ and $g\in L^1(L)$:
\begin{equation*}
\int_{L/M}f(lM)\, d(lM) =\int_{L}f(lM)\, dl=
\int_L (f\circ\pi)(l)\, dl
\end{equation*}
and
\begin{equation*}
\int_{L} g(l)\, dl=\int_{L/M}\int_{M}g(lm)\; dm\, d(lM)\, ,
\end{equation*}
where $\pi: L\to L/M$ is the canonical projection $l\mapsto lM$.

Let $B:=\exp(\mathfrak{b})$ and $B_1:=\exp (\fb_1)=B\cap U^\prime$.
Then $U=KBK=TKB_1K$.
In particular, denoting by $x_{0}$ the point $\{K\}\in U/K$,
then $KB\cdot x_{0}=T(KB_1)\cdot x_0=U/K$.
Set $M=Z_{K}(B)$ and define $\Psi:K/M\times
B\rightarrow U/K$ by $\Psi(kM,b):=kb\cdot x_{0}$.
Then $\Psi$ is smooth and surjective. Furthermore,
\begin{equation} \label{eq:delta}
|\det(d\Psi_{(kM,\exp(H))})|=
\prod_{\alpha\in\Delta^{+}}|\sin\alpha(H)|^{m_{\alpha}}=:\delta(\exp(H)).
\end{equation}
This proves the following integration formula
(cf. e.g. \cite{He84}, Theorem 5.10.)

\begin{lemma}
There exists a constant $c>0$ such that for all $f\in C(U/K)$ we have
\begin{equation*}
\int_{U/K}f(uK)\, d(uK)=c\int_{K/M}\int_{B}f(kb\cdot x_{0})~\delta(b)\;
db \, d(kM)\, .
\end{equation*}
\end{lemma}

\section{Spherical functions and spherical representations}
\label{section:spherical}

\noindent
In this section we recall some necessary facts about
spherical functions and spherical representations.
We refer to \cite{He84}, Chapter V, as standard reference.
Our main result in this section is Theorem \ref{th-extension}. It states
that, if $m_\alpha$ is even
for all $\alpha\in \Delta$, then there exists a differential
operator $D$ on $B$ with analytic coefficients
and a rational function $Q$  on
$\fa_\C^*=\fb_\C^*$  such
that
$$\delta (b) \psi_\mu (b)= Q(\mu)D\left(\sum_{w\in W}
b^{w(\mu+\rho)}\right)\, .$$
Here $\psi_\mu$ is the spherical function corresponding to the spherical
representation with highest weight $\mu$,  the function $\delta$
is as in (\ref{eq:delta}), and 
\begin{equation} \label{eq:rho}
\rho=\frac{1}{2}\sum_{\alpha\in \Delta^+}
m_\alpha \alpha\in \fa^*.
\end{equation}
The differential operator $D$ will be constructed explicitly 
from Opdam's shift operator and in addition the rational function $Q$ will
be explicitly determined. It is a holomorphic extension (and a
trivial extension to $T$) of the differential
shift operator constructed in \cite{OP04} for the noncompact symmetric 
case with even multiplicities.

\subsection{Spherical representations}
Let $(\pi,V)$ be an irreducible unitary representation of $U$. Let
\begin{equation*}
V^{K}:=\{v\in V\mid\forall k \in K~:~\pi(k)v=v\}~.
\end{equation*}
We say that $\pi$ is \textit{spherical} if $V^{K}\not =\left\{ 0\right\}$.
In this case, then $\dim V^K=1$.
We denote by $\widehat{U}$ the set of equivalence classes of
irreducible unitary
representations of $U$ and by $\widehat{U_{K}}$ the subset of equivalence
classes of irreducible $K$-spherical representations. We shall use the same
notation
for a given unitary representation and for the corresponding
equivalence class in $\widehat{U}$.

If $\lambda \in \fb_\C^*$ and $b=\exp (H)\in B$, then we
write $b^\lambda =e^{\lambda (H)}$ provided this is well defined.
The same notation will be adopted for elements in the complexification of $B$.
Let $\pi\in \widehat{U}$. As $T$ is central in $U$, it follows
by Schur's Lemma that $\pi$ is of the form
\begin{equation*}
\pi(tu)=t^\lambda \pi^{\prime}(u)~,\qquad t\in T,\,u\in U^{\prime}\,,
\end{equation*}
where $\lambda$ is some element of $i\fb^*$ and $\pi^\prime=\pi|_{U^\prime}$.

If $\fz\not= \{0\}$, then we let $\Gamma_0 :=\{X\in \fz\mid\exp (X)=e\}$, as
before. Then
$$i\Gamma_0^*:=\{\lambda \in i\fz^*\mid
\forall H\in \Gamma_0\, :\, \lambda (H)\in 2\pi i \Z\}\simeq \widehat{T}\,,$$
where the isomorphism is given by
$\lambda \mapsto \chi_\lambda$ and $\chi_\lambda (t):=t^{\lambda}$.
Note that, if we do not assume  $T\cap U^\prime=\{e\}$, then
we have to impose the additional condition that $\pi^\prime (t)=
t^\lambda \id$ for all $t\in T\cap U^\prime$.

Let $\mathfrak{c}$ be a Cartan subalgebra of $\mathfrak{u}$
containing $\fb$. Set $\fc_1:=\fc\cap \fu^\prime$.  We say that $\mu\in
i\fc^{\ast}$ is an \textit{extremal weight} of an irreducible
representation $\pi$ of $U$ if $\mu$ is the
highest weight of $\pi$ with respect to some ordering in
$i\fc^\ast$.
We fix an ordering on $i\fz^*$, then we extend it
to $i\fb^*$ by using the lexicographic ordering on
$\fa_1^*$, and we finally extend it to an ordering in 
$i\mathfrak{c}^{\ast}$. 
If $\pi$ is an
irreducible representation of $U$, then $\mu(\pi)\in i \fc^*$
denotes the highest weight of $\pi$ with respect to this ordering. Similarly,
if $\sigma\in \widehat{U^\prime}$, then $\mu (\sigma )\in i\fc^*_1$ denotes
the highest weight of $\sigma$.
Notice that, in this notation, we have
$$\mu (\pi) |_{\fc_1}=\mu (\pi |_{U^\prime})\, .$$
For
$\lambda\in \mathfrak{c}^*_\C$ and $\ga\in \Delta$ let
\begin{equation}\label{lambdaa}
\lambda_\alpha =\frac{\inner{\lambda}{\ga}}{\inner{\ga}{\ga}}
=\frac{1}{2}\lambda(H_\alpha)\, .
\end{equation}
Denote by $\gL_K^+=\gL_K^+(U)$ the set of highest weights of spherical 
representations
of $U$. Then we have $\gL_K^+(U)=i\gG^*_0\oplus \gL_K^+(U^\prime)$.
Here we employ the notation $\oplus$ to indicate that an element 
of $\gL_K^+(U)$ can be written in a unique way as a sum of an element
of $\gG^*_0$ and an element of $\gL_K^+(U^\prime)$.
If $\mu \in \gL_K^+(U)$, then $\pi_\mu$ denotes the corresponding
spherical representation and $w_\mu$ a $K$-invariant vector in the 
space of $\pi_\mu$ satisfying $\|w_\mu \|=1$.

\begin{theorem}\label{th-spher}
Let $(\pi,V)$ be an irreducible representation of $U$. Then
the following holds:

\begin{enumerate}
\item If $\pi$ is spherical then $\mu (\pi)\in i\gG^*_0 \oplus \fa^{\ast}_1$
and
\begin{equation}\label{eq-spher}
\frac{\inner{\mu (\pi)}{\alpha}}{\inner{\alpha}{\alpha}}=:
\mu (\pi)_\alpha \in\N_0 
\end{equation}
for all $\alpha\in\Delta^+$. Here $\N_0=\{0,1,2,\dots\}$.

\item Let $\mu\in i\gG^*_0 \oplus \mathfrak{a}^{\ast}_1$ so that
$\mu_{\ga} \in \Z$ for all $\ga\in \gD$. If $U^\prime$ is
simply connected, then there exists a unique
spherical representation $\pi$ with extremal
weight $\mu$.

\item If $U^\prime$ is simply connected
then $\Lambda^+_K(U)= i\gG^*_0\oplus \{\mu \in \fa^*_1
\mid \forall \alpha\in\gD^+\, :\, \mu_\alpha\in \N_0\}\, .$
\end{enumerate}
\end{theorem}

\subsection{Spherical functions}
Recall the following definition.

\begin{definition}
Let $G$ be a locally compact Hausdorff topological group 
and $K\subset G$ a compact subgroup. 
A continuous function $\varphi: G\to\mathbb{C}$ is said to be 
\emph{spherical} if $\varphi$ is $K$-bi-invariant, is not identically $0$,
and satisfies the identity
\begin{equation*}
\int_{K} \varphi(xky)\, dk =\varphi(x)\varphi(y)
\end{equation*}
for all $x,y\in G$.
\end{definition}

For $\mu\in \gL^+_K(U)$ define
$\psi_{\mu}:U\rightarrow\mathbb{C}$ by%
\begin{equation}\label{eq-coeff}
\psi_{\mu}(u )=(\pi_\mu (u)w_\mu ,w_\mu)\, ,
\end{equation}
where $(\cdot,\cdot)$ denotes the inner product in the space of 
$\pi_\mu$ for which this representation is unitary.
Then $\psi_\mu$ is a spherical function on $U$ and
every spherical function on $U$ is of the form $\psi_\mu$
for some $\mu\in \gL^+_K(U)$. 
Notice that, with $\lambda:=\mu|_{\fz}\in \gG^*_0$
and $\mu^\prime:=\mu|_{\fb_1}\in \gL^+_K(U^\prime)$, we have
\begin{equation}\label{eq-sptb}
\psi_\mu(tu')=t^{\lambda} (\pi_{\mu^\prime}(u')w_\mu ,w_\mu)
=t^\lambda \psi_{\mu^\prime} (u'), \qquad  t\in T\,,\; u'\in U'.
\end{equation}

Since $\pi_\mu$ is unitary, (\ref{eq-coeff}) implies the following lemma. 

\begin{lemma}\label{le-ginverse} 
Let $\mu\in \gL^+_K(U)$. Then
$$\overline{\psi_\mu (u)}=\psi_\mu (u^{-1})$$
for all $u \in U$.
\end{lemma}

Let $\iota: U\to U_\C$ be the universal complexification of
$U$. Hence, if $L$ is a complex Lie group
and $\varphi : U\to L$ is a Lie group homomorphism, then there
exists a holomorphic homomorphism $\varphi_\C : U_\C\to L$ such
that $\varphi_\C\circ \iota = \varphi$. As $U$ is compact, it
follows that there exists a faithful representation
$\pi : U\to \GL (n,\C)$ for some $n$. Applying the above to $\pi$, 
we conclude that $\iota$ has to be injective. 
We can therefore assume that
$U$ is a subgroup of $U_\C$. Since $U$ is compact, it follows that
$U$ is closed in $U_\C$.

\begin{lemma}\label{le-holext} 
Let $\mu\in \Lambda^+_K(U)$. Then the
spherical function $\psi_\mu$ extends to a holomorphic function
on $U_\C$. The extension is given by
$$\psi_\mu (g)=((\pi_\mu)_\C (g)w_\mu ,w_\mu )\, .$$
\end{lemma}

Let $G$ to be the analytic subgroup of
$U_\C$ with the Lie algebra $\fg :=\fk\oplus i\fq$. Then
$G$ is closed in $U_\C$ and $K\subset G$.
We set $\fp:=i\fq$ and notice that
$\fa=i\fb$ is a maximal abelian subspace of $\fp$.
Denote  by $\tau_\C$ the holomorphic
extension of $\tau$ to $U_\C$ and set $\theta=\tau_\C|_G$.
Then $\theta$ is a Cartan involution on $G$. We have
$K=G^\theta$. The symmetric
space $G/K$ is called a \textit{noncompact dual} of $U/K$.
We set $A=\exp (\fa)$ and $A_1=\exp (\fa_1)$.
Finally we set $T_\C=\exp (\fz_\C)$ and $T_\R=T_\C\cap G=\exp (i\fz) 
\subseteq A$.

Let us recall the standard notations and definition for
the Iwasawa decomposition of $G$.
Let
$$\fn=\oplus_{\alpha\in\gD^+} \fg^\alpha\, ,$$
where, as usual,
$\fg^\ga=\{X\in \fg\mid \forall  H\in \fa : [H,X]=\alpha (H)X\}$.  
Then the Iwasawa map
\begin{equation}\label{eq-iwasawa}
K\times A\times N\ni (k,a,n)\mapsto kan\in G
\end{equation}
is a diffeomorphism. For $x\in G$ define $(k(x),a(x),n(x))\in K\times A\times
N$
by the inverse of the map in (\ref{eq-iwasawa}).
We normalize the Haar measure $dn$ on $N$ so that $\int_{N} a (\theta
(n))^{-2\rho}\, dn=1$.
Here $\rho$ is as in (\ref{eq:rho}).
Moreover, we normalize the Haar measure $da$ on $A$ (and similarly $dt$
on $T_\R$) so that the \textit{Fourier transform} on $C^\infty_c(A)$,
defined by 
\begin{equation*}
\widehat{f}(\lambda) =\mathcal{F}_{A}(f)(\lambda)
:=\int_{A} f(a)a^{-\lambda}\, da\,, \qquad f \in C^\infty_c(A)\,,\;
\lambda \in \fa^*_\C\,,
\end{equation*}
has inverse
$$f(a)=\int_{i\fa^*} \widehat{f}(\lambda )a^\lambda\, d\lambda\, .$$
Finally we normalize the Haar measure $dg$ on $G$ so that the equality
$$\int_G f(g)\, dg =\int_K\int_{A} \int_N f(kan)a^{2\rho}\;
dn\,da\,dk$$
holds for all $f\in C_c(G)$.

For $\lambda\in \fb_\C^*=\fa_\C^*$ let
\begin{equation}\label{eq-sphericalfct}
\varphi_\lambda (g)=\int_K a(g^{-1}k)^{-\lambda -\rho}
\, dk=\int_K a (gk)^{\lambda -\rho}\, dk\,
\end{equation}
be the corresponding spherical function on $G$. Let 
$G^\prime$ be the analytic subgroup of $G$ corresponding to the
Lie algebra $\fg^\prime:=[\fg,\fg]$.
Notice that if $g=th$ with $t\in T_\R$ and $h\in G^\prime$, then
$$\varphi_\lambda (g)=t^{\lambda_1}\varphi_{\lambda_2}(h)$$
where 
$\lambda_1$ is the restriction of $\lambda$ to $\fz_\C$, $\lambda_2$
is the restriction of $\lambda$ to $\fa_{1\C}=\fa_\C\cap [\fu_\C ,\fg_\C]$,
and $\varphi_{\lambda_2}$ is the spherical function
on the semisimple Lie group $G^\prime$ with spectral
parameter $\lambda_2$.
Recall that $\varphi_\lambda =\varphi_\mu$ if and only
if there exists $w\in W$ such that $\lambda =w\mu$.
We also recall the following well known fact.

\begin{lemma}
\label{le:sphercnc} 
Let $\mu\in\Lambda_{K}^{+}(U)$ and
let $\psi_{\mu}$ denote the holomorphic extension to $U_\C$ of
the spherical function $\psi_\mu$ on $U$.
Then
\begin{equation*}
\psi_{\mu}|_G=\varphi_{\mu+\rho}~.
\end{equation*}
\end{lemma}
\begin{proof} (Cf. e.g. \cite{He84}, pp. 540-541.)
Fix a highest weight vector $u$ for $\pi_\mu$ such that
$w_\mu=\int_{K}\pi_\mu (k)u\, dk$. In particular
$(u, w_\mu)=1$ and for $b\in B$ we have
\begin{equation}\label{eq-int}
\psi_{\mu}(b) =(\pi_\mu(b)w_\mu, w_\mu)
=\int_{K}(\pi_\mu(bk)u , w_\mu) \, dk.
\end{equation}
As $K$ is compact, it follows that (\ref{eq-int}) remains valid for the
holomorphic extension of $\psi_{\mu}$. 
In particular it is valid for $b\in A$.
Thus, as $(u,w_\mu)=1$,
\begin{equation*}
(\pi_\mu(bk)u,w_\mu)=a(bk)^{\mu}(u , w_\mu)=a(bk)^{(\mu+\rho)-\rho}
\end{equation*}
and hence
\begin{equation*}
\psi_{\mu}(b)=\int_{K}a (bk)^{(\mu+\rho)-\rho}\;dk
=\varphi_{\mu+\rho}(b)~.
\end{equation*}
\end{proof}

\subsection{The dimension function $d(\mu)$ and the $\bc$-function}
Set $\overline{N}:=\theta(N)$ and normalize the Haar measure 
$d\overline{n}$ on $\overline{N}$ so that 
$\int_{\overline{N}} a (\overline{n})^{-2\rho}\, d\overline{n}=1$.
If $\mathrm{Re} \lambda_\alpha > 0$ for all $\alpha \in \gD^+$, 
then the Harish-Chandra
$\bc$-function for $G/K$ is
given on $\fa_\C^*$ by
\begin{equation}\label{eq-cfunction}
\bc(\lambda )=\int_{\bar{N}}
a(\bar{n})^{-\lambda -\rho}\, d\bar{n}\,,
\qquad \lambda \in \fa^*_\C.
\end{equation}
Observe that $\bc(\rho)=1$. By the
Gindikin-Karpelevic formula we have
\begin{equation*}
\bc(\lambda)=
c_{0}\prod_{\alpha\in\Delta^{+}}\bc_{\alpha}(\lambda_{\alpha})\,,
\end{equation*}
where $\bc_{\alpha}(\lambda_{\alpha})$ corresponds to a rank-one
$\bc$-function, i.e.
\begin{equation}\label{cfunct}
\bc_\alpha (\lambda_\alpha)=
\frac{2^{-\lambda_{\alpha}}
\Gamma\left(\lambda_{\alpha}\right) }
{\Gamma(\frac{1}{2}(\frac{1}{2}m_{\alpha}+1+\lambda_{\alpha}))
\Gamma(\frac{1}{2}(\frac{1}{2}m_{\alpha}+m_{2\alpha}+\lambda_{\alpha}))}\,,
\end{equation}
and the constant $c_{0}$ is determined by $\bc(\rho)=1$.
In particular, this formula gives the meromorphic extension
of $\bc$ to all
of $\fa_\C^*$. If
$m_{2\alpha}=0$ for all $\alpha$ then (\ref{cfunct}) simplifies to
\begin{equation*}
\bc(\lambda)=c_{1}\prod_{\alpha\in\Delta^{+}}
\frac{\Gamma(\lambda_{\alpha})}{\Gamma(\lambda_{\alpha}+m_{\alpha}/2)}\,.
\end{equation*}

If $m_\alpha$ is even for all $\alpha\in \Delta$, then $m_{2\alpha}=0$
and, by the classification or by \cite{OP04}, Appendix C, there exists 
$m\in \N$ such that
$m_\alpha=2m$ for all $\alpha\in \gD$. The
relation $\Gamma(z+1)=z\Gamma(z)$ implies then the following lemma.

\begin{lemma}\label{c-function}
Suppose that $m_{\alpha}$ is even for all $\alpha\in\Delta$. Let
$2m\in 2\N$ be the resulting common value of $m_\alpha$ for all
$\alpha\in\gD$.
Then
\begin{equation*}
\frac{1}{c(\lambda)}=C\prod_{\alpha\in\Delta^{+}}
\prod_{k=0}^{m-1}(\lambda_{\alpha}+k)
\end{equation*}
where the constant $C$ is given by
$$C=\prod_{\alpha\in\Delta^+}\prod_{k=0}^{m-1}\frac{1}{\rho_\alpha+k}\, .$$
\end{lemma}

The dimension $d(\mu)$ of the spherical representation 
$\pi_\mu$ can be expressed as a limit of ratios of $\bc$-functions by means of 
Vretrare's formula, cf.\ \cite{He94}, Theorem 9.10, p.\ 337. 
In the even multiplicity case this formula simplifies because the limit 
involved can be computed as the quotient of the limits of the 
$\bc$-functions appearing in the numerator and in the denominator of 
the formula.

\begin{lemma}\label{le-dim}
Assume that for all $\alpha\in\gD$ the multiplicities $m_\alpha$ 
are even, and let $2m$ be their resulting common value.
Then the following properties hold:
\begin{enumerate}
\renewcommand{\labelenumi}{(\theenumi)}
\item $\rho_\alpha=m$ for every simple root $\alpha$\,;
\item $\rho_\alpha\in\Z$ for every $\alpha \in\Delta$\,;
\item $\rho_\alpha\geq m$ for every $\alpha \in\Delta^+$\,;
\item For all $\mu\in \gL_K^+(U)$ we have
$$
d(\mu)=\frac{\bc(-\rho)}{\bc(\mu+\rho)\bc(-(\mu+\rho))}\,.
$$
\end{enumerate} 
\end{lemma}
\begin{proof}
If $\alpha$ is a simple root in a reduced root system
$\Delta$, then $\rho_\alpha=m_\alpha/2$. Indeed, 
$\rho-(m_\alpha/2)\alpha$ is fixed by the reflection $s_\alpha$. This 
proves (1). All the remaining statements follow easily from the
first and from Vretrare's formula.
\end{proof}

\begin{theorem}\label{th-cfunctdim}
Assume that all $m_\alpha$ are even for all $\alpha\in\gD$,
and let $2m$ be their resulting common value.
Then the dimension function $d$ extends as to a polynomial function on 
$\fa^*_\C$ given by 
\begin{equation*}
d(\lambda )
=\prod_{\alpha\in \Delta^{+}}\prod_{k=0}^{m-1}
\frac{k^{2}-(\lambda +\rho)_{\alpha}^{2}}{k^2-\rho_\alpha^2}\, .
\end{equation*}
\end{theorem}
\begin{proof} This follows from Lemmas \ref{c-function} and \ref{le-dim}.
\end{proof}

\subsection{The differential shift operator}
Suppose that all multiplicities $m_\alpha$ are even. Then the
function $\delta(b)$ of (\ref{eq:delta}) extends as a 
$W$-invariant holomorphic function on $B_\C=AB$. (See Lemma
1.2 in \cite{OP04}.) The following theorem, which is a slight extension
of Theorem 5.1(c) of \cite{OP04}, provides an explicit formula on $A$
for the spherical functions on a symmetric space $G/K$. This theorem is 
our starting point for investigating the holomorphic extension of the 
spherical functions on $G/K$ to its complexification $U_\C/K_C$. 
By restriction, we shall then deduce explicit formulas for the spherical 
functions on the compact dual symmetric space $U/K$. 
Because of our context, 
we shall only consider the case for which
$m_\alpha$ is even for each 
$\alpha \in \gD$.  Recall that since we also assume that $\fu^\prime$
is simple, this means that the $m_\alpha$ have a common value $2m\in 2\N$.

\begin{theorem}\label{th-extension} 
Assume that all multiplicities $m_{\alpha}$ are even.
Then there exists a $W$-invariant differential operator $D$ on $A$  
with analytic
coefficients, such that for all $\lambda\in {\fb}_{\mathbb{C}}^{\ast}=
{\fa}_{\mathbb{C}}^{\ast}$ and all $a \in A$ we have
\begin{equation}\label{eq-extension}
\delta(b)\varphi_{\lambda}(b)
 =\frac{1}{d(\lambda -\rho )}\, D\left( \sum_{w\in W}b^{w\lambda}\right).
\end{equation}
The right hand side of (\ref{eq-extension}) is holomorphic in $\lambda$.
\end{theorem}
\begin{proof} 
If $G$ (and hence $U$) is semisimple, then this
follows from Lemma \ref{le-dim}, Theorem \ref{th-cfunctdim},
and \cite{OP04}, Theorem 5.1(c). For the
general case, let $D^\prime$ be the differential operator from
\cite{OP04} on $A_1$. By assumption, we have $A\simeq T_\R\times A_1$. 
For $t\in T_\R$ and $a\in A_1$ we define the operator $D$ by
$D(f)(ta):=
D_a^\prime f(ta)/\bc(-\rho)$,
where the subscript 
$a$ indicates differentiation with respect
to the variable $a$.
\end{proof}

We remark that
the operator $D'$ occurring in the proof of Theorem \ref{th-extension}
is of the form $\delta(a) \widetilde{D}$, where $\widetilde{D}$ is
Opdam's shift operator of shift $2m$. Multiplication by $\delta(a)$
cancels the singularities of the coefficients of $\widetilde{D}$.
By construction, $D'$ can be considered as differential
operator on $B_1$ with holomorphic coefficients on $(B_1)_\C=A_1B_1$.
Consequently the operator $D$ itself can be considered as differential
operator on $B$ with holomorphic coefficients on $B_\C=AB$.

The following corollary, which allows us to holomorphically extend the 
right hand side of (\ref{eq-extension}), will also 
play a crucial role in the proof of the local Paley-Wiener theorem.

\begin{corollary} \label{co-zero}
Suppose that $\lambda\in\mathfrak{b}_{\mathbb{C}}^{\ast}$ is
such that $\lambda_{\alpha}\in\pm \{0,1,\ldots,m-1\}$ for some 
$\alpha\in\Delta^{+}$. Then $D\left( \sum_{w\in W}b^{w\lambda}\right)=0$
for all $b \in B_\C:=AB$. 
\end{corollary}

\begin{proof}
The last statement in Theorem \ref{th-extension} ensures that 
$D\left( \sum_{w\in W}b^{w\lambda}\right)=0$ for all $b \in A$.
Since $D\left( \sum_{w\in W}b^{w\lambda}\right)$ is holomorphic
(possibly multivalued) on $B_\C$ and 
and it vanishes on $A$, it must be identically zero on $B_\C$.
\end{proof}

Because of Corollary \ref{co-zero}, the only obstruction to 
the holomorphic extension to $B_\C$ of the right hand side of 
(\ref{eq-extension}) is the fact that the functions $b^{w\lambda}$
might be multivalued. This obstruction is solved by choosing
a domain where the exponential function is a diffeomorphism. 

Let $0\in \cV\subset \fb_\C$ and
$e\in \cU\subset B_\C$ be open, connected and such that $\exp :\cV\to
\cU$ is an analytic diffeomorphism. We will assume furthermore that
$\cV\cap \fb$ is open and connected, and
that $\cV$ (and hence also $\cU$) is
$W$-invariant and contains $\fa$. Then, by Theorem \ref{th-extension},
and Corollary \ref{co-zero}, $\varphi_\lambda$ has an analytic 
extension to $\cU$.

\begin{theorem}\label{th-anaext}  The function
$$\gL^+_K(U)\times \cU\cap B \ni (\mu,b)\mapsto
\psi_\mu (b)\in \C$$
has a holomorphic extension to  $\fb_\C^\ast\times \cU$ given by:
\begin{eqnarray*}
\psi_{\lambda}(b)  & =&\varphi_{\lambda+\rho}(b)\\
&=& \frac{1}{d(\lambda )\delta (b)}\,D\left(\sum_{w\in W}b^{w(\lambda +\rho)}
\right)\\
&=&\frac{1}{d(\lambda )\delta (b)}\sum_{w\in W}Db^{w(\lambda+\rho)}\\
&=& \delta(b)^{-1}\left(
\prod_{\alpha\in\Delta^{+}}\prod_{k=0}^{m-1}
\frac{k^2-\rho_\alpha^2}{k^{2}-(\lambda +\rho)_{\alpha}^{2}}\right)
 ~D\left( \sum_{w\in W}b^{w(\lambda + \rho)}\right)\, .
\end{eqnarray*}
Furthermore the analytic continuation satisfies
$$\psi_\lambda = \psi_{w(\lambda +\rho)-\rho}$$
for all $w\in W$ and $\lambda \in \fb_\C^*$.
\end{theorem}

\begin{proof} The first claim follows from Lemma \ref{le:sphercnc},  
Theorem \ref{th-extension} and Corollary \ref{co-zero}, as the right 
hand side of (\ref{eq-extension})
extends to a analytic function on $\fb_\C^*\times \cU$. 
The statement $\psi_\lambda =\psi_{w(\lambda +\rho)-\rho}$
follows from the Weyl group invariance of 
$\lambda \mapsto \varphi_\lambda (b)$.
\end{proof}

As a corollary, we obtain the following explicit formulas for the 
spherical functions on Riemannian symmetric spaces of compact type
and even multiplicities.

\begin{corollary}\label{cor:psi}
Let $\mu\in\Lambda_{K}^{+}(U)$. Suppose that all multiplicities $m_\alpha$ 
are even.
Then the following holds on $B$:
\begin{equation} \label{eq:spherc}
\delta(b)\psi_{\mu}(b)=\frac{1}{d(\mu)}D\left(\sum_{w\in W}b^{w(\mu+\rho
)}\right)=\frac{1}{d(\mu)}\sum_{w\in W}Db^{w(\mu+\rho)}~.
\end{equation}
\end{corollary}

\begin{proof} This follows from the fact that $w(\mu + \rho) \in
 \Lambda_{K}^{+}(U)$ and that for all $\nu \in \Lambda_{K}^{+}(U)$
the function $b^\nu$ is single valued and holomorphic on $B_\C$.
\end{proof}

\subsection{The classification}
\label{subsection:classification}
We finish this section by giving the classification -- up to coverings -- of
the symmetric spaces $U/K$ with $U$ semisimple and $U/K$ irreducible.
Here $K$ stands for an arbitrary
connected, compact and simple Lie group. We list
also the noncompact Riemannian dual $G/K$ as well as
$r:=\mathrm{rank}(U/K)=\dim(\mathfrak{b})$, the
multiplicities $m_\alpha$, and the dimension $d$
of $U/K$. In all these cases the multiplicities $m_\alpha$ are
constant.

\medskip

\begin{table}[h]
\label{table:UmodK}
\setlength{\extrarowheight}{2pt}
\begin{tabular}{|c|c|c|c|c|}
\hline
\multicolumn{5}{|c|}{
\rule[-2mm]{0mm}{7mm}\textsf{Table I}}\cr
\hline\hline
$U/K$ &$G/K$ & $m_{\alpha}$ & rank $r$ & dimension $d$\\[.2em]
\hline
$K\times K/K\simeq K$ & $K_{\mathbb{C}}/K$& 2 & $r$ & $d$ \\[.2em]
\hline
$\mathrm{SU}(2n)/\mathrm{Sp}(n)$ &$\mathrm{SU}^*(2n)/\mathrm{Sp}(n)$  
& 4 & $n-1$  &$(n-1)(2n+1)$ \\[.2em]
\hline
$\mathrm{E}_{6(-78)}/\mathrm{F}_{4}$ & $\mathrm{E}_{6(-26)}/\mathrm{F}_4$& 8 &
2 & 26 \\[.2em]
\hline
$\stackrel{\ }{\displaystyle\frac{\mathrm{SO}(2(n+1))}{\mathrm{SO} _0 (2n+1)}}
\simeq S^{2n+1}$
 &$\stackrel{\
 }{\displaystyle\frac{\mathrm{SO}_0(2n+1,1)}{\mathrm{SO}(2n+1)}}$
& $2n$ & 1 &
$2n+1$ \\[.3em]
\hline
\end{tabular}
\end{table}

The first line of Table I corresponds to the complex case, in which the
Lie algebra $\fg$ admits a complex structure. The classification of
compact real forms $K$ for simple complex Lie algebras $\fg=\fk_\C$
is recorded in Table II (see \cite{He84}, p.\ 516).

\begin{table}[h]
\label{table:KwithGcomplex}
\setlength{\extrarowheight}{4pt}
\begin{tabular}{|c|c|c|c|}
\hline
\multicolumn{4}{|c|}{
\rule[-2mm]{0mm}{7mm}\textsf{Table II}}\cr
\hline\hline
$K$ &$d:=\dim K$ & rank $r$ & $\fg$\\[.2em]
\hline
$\mathrm{SU}(n)$ & $(n-1)(n+1)$ & $n-1$ & $\fa_{n-1}$\\[.2em]
\hline
$\mathrm{SO}(2n)$ & $n(2n-1)$ & $n$ &$\fd_{n}$ \\[.2em]
\hline
$\mathrm{SO}(2n+1)$ & $n(2n+1)$ & $n$ &$\fb_{n}$ \\[.2em]
\hline
$\mathrm{Sp}(n)$ & $n(2n+1)$ & $n$ &$\fc_{n}$ \\[.2em]
\hline
$\mathrm{E}_6$ & 78 & 6 & $\fe_6$ \\[.2em]
\hline
$\mathrm{E}_7$ & 133 & 7 & $\fe_7$ \\[.2em]
\hline
$\mathrm{E}_8$ & 248 & 8 & $\fe_8$ \\[.2em]
\hline
$\mathrm{F}_4$ & 52 & 4 & $\ff_4$ \\[.2em]
\hline
$\mathrm{G}_2$ & 14 & 2 & $\fg_2$ \\[.2em]
\hline
\end{tabular}
\end{table}

\section{The local Paley-Wiener Theorem for
compact symmetric spaces with even multiplicity}
\label{section:localPW}

\noindent
In this section we introduce the spherical Fourier transform of
$K$-invariant functions on the compact symmetric
spaces $U/K$. We then use the results from the last section,
in particular Theorem \ref{th-anaext}, to show that the
Fourier transform, which in the beginning is only defined on
a discrete set, extends holomorphically to $\fb_\C^*$ as long as
the $K$-invariant function has sufficiently small support. We
then define the Paley-Wiener space on $\fb_\C^*$ and
prove the local Paley-Wiener theorem.  This theorem generalizes 
the results obtained by Gonzalez in \cite{Go00} for the case
$U=K\times K$, where $K$ is a connected, compact,
simple Lie group. Notice that in this case
$U/K\cong K$.

Let $\|\cdot \|$ be the norm on $\fu$ with respect to
the $U$-invariant inner product constructed in Section 1, and let
$d$ be the associated Riemannian distance function on $U/K$.
For $R>0$ let $B_R:=\{X\in \mathfrak{q}\mid \|X\|\le R\}$ and 
$D_R:=\{x \in U/K\mid d(x,x_0) \leq R\}$
denote the corresponding balls of radius $R$ with center $0$ and 
$x_0$ respectively. We suppose that $R$ is chosen so that the map
$\Exp: X \to \exp X \cdot x_0$ is a diffeomorphism of $B_R$ onto $D_R$.
Finally we define
\begin{equation}\label{def-CR}
C^\infty_R(U/K)^K:=\{f\in C^\infty (U/K)^K \mid \Supp (f)\subseteq D_R\}\,.
\end{equation}
Here and in the following the superscript $K$ denotes $K$-invariance.

Note that $\pi : B\to B\cdot x_0$ is a finite covering. We
will  identify $B$ with the image $B \cdot x_0$. This
is allowed, as we will only be considering
$K$-invariant functions. Then, for every $f \in C^\infty_R(U/K)^K$
we have $\Supp f \subseteq D_R$ if and only if $\Supp(f|_B) \subseteq D_R$.

\subsection{The spherical Fourier transform on $U/K$}
For $f\in L^{2}(U/K)^{K}$ define
$\widehat{f}:\Lambda_{K}^{+}(U)\rightarrow \mathbb{C}$ by
\begin{equation}\label{def-Four}
\begin{array}{rl}
\mathcal{F}(f)(\mu) = \widehat{f}(\mu)
&\!\!\!:=(f,\psi_\mu ) \\
&\!\!\!= \displaystyle\int f(u)\psi_{\mu}(u^{-1})\, du  \\
&\!\!\!= c\displaystyle\int_{B}f(b)\psi_{\mu}(b^{-1})\delta(b)~db
\end{array}
\end{equation}
Here we have used the equality 
$\overline{\psi_\mu (u)}=\psi_\mu (u^{-1})$
 from Lemma \ref{le-ginverse}, and $c$ is a suitable positive constant
depending on the fixed normalization of measures.
We call $\widehat{f}$ the \textit{spherical Fourier transform} 
of $f$ and the map $\mathcal{F}$ the \textit{spherical 
Fourier transform}.  It
is well known that $\cF$ is a unitary isomorphism of
$L^2(U/K)^K$ onto $\ell^2(\gL^+_K(U))$. The inversion formula is stated
in the next theorem. In the following we shall often consider the 
$K$-bi-invariant functions $\psi_\lambda$ on $U$ as $K$-invariant functions
on $U/K$.

\begin{theorem}\label{th-finv}
Let $f\in L^{2}(U/K)^{K}$. Then
\begin{equation*}
f=\sum_{\mu\in\Lambda_{K}^{+}}d(\mu )\widehat{f}(\mu)\psi_{\mu}~,
\end{equation*}
where the sum is taken in $L^{2}(U/K)^{K}$. If $f\in C^{\infty}(U/K)^{K}$
then the above sum converges in the $C^\infty$-topology.
\end{theorem}

Let $E$ be a differential operator on $B$. Then the \textit{formal adjoint}
operator $E^{\ast}$ is defined by the relation 
\begin{equation*}
\int_B f(b)Eg(b)\, db =\int_B E^\ast f(b)g (b)\, db
\end{equation*}
for all $f,g\in C^\infty(B)$. In this section $D$ denotes
the differential operator
from Theorem \ref{th-anaext}. Hence
\begin{equation}
\delta (b)\psi_\mu (b)=d(\mu )^{-1} \sum_{w\in W}D b^{w(\mu +\rho)}\, .
\end{equation}
If $\lambda = \lambda_R+ i\lambda_I\in \fb_\C^*$ with 
$\lambda_R,\lambda_I\in \fb^*$,
then we set 
\begin{equation} \label{eq:ReIm}
\re\lambda:=\lambda_R\qquad \text{and} \qquad \im\lambda:= \lambda_I\,.
\end{equation}
\begin{definition} 
A holomorphic function $f:\fb_\C^*\to \C$ is said to
be of \emph{exponential type $R>0$} 
if for each $N\in \N_0$ there exists a constant $C_N>0$ such that 
$$|f(\lambda )|\le C_N(1+\|\lambda \|)^{-N}e^{R\|\re \lambda \|}$$
for all $\lambda\in\fb_\C^*$.
\end{definition}

\begin{lemma}\label{le-exp} 
Let $p:\fb_\C^*\to \C$ be a polynomial function. 
Assume that $F:\fb_\C^*\to \C$ is an entire function 
so that $\lambda \mapsto p(\lambda )F(\lambda )$ is of 
exponential type $R$. 
Then $F$ itself is of exponential type $R$.
\end{lemma}
\begin{proof} 
Let $n=\dim_\C \fb_\C^*$. Choose a basis $\lambda_1,\ldots ,\lambda_n$
of $\fb_\C^*$ and identify $\fb_\C^*$ with $\C^n$ by
$\sum z_j\lambda_j\mapsto (i z_1,\ldots ,iz_n)$. 
The fact that $F$ is of exponential type follows then from 
\cite{He94}, Ch.\ III, Theorem 5.13. The type $R$, which is not explicitly
computed in this reference, can be easily deduced from formula (62) in that
proof. Indeed one obtains the estimate
$$|F(z)|\le C \sup_{\|w\|\le 1}|p(z+w)F(z+w)|\,, $$
which shows that $F$ and $p F$ have the same exponential type.
\end{proof}

\begin{theorem}\label{le-pw2}
Suppose that all multiplicities are even.
Let $\mu\in\Lambda_{K}^{+}(U)$ and $f\in C^{\infty}(U/K)^{K}$. Then
\begin{equation*}
\widehat{f}(\mu)  =\frac{c}{d(\mu)} |W|\int [D^{\ast}f(b)]\, b^{-w(\mu
+\rho)}\, db \,,
\end{equation*}
where $|W|$ denotes the cardinality of $W$.

Assume that $R>0$ is chosen so that $\Exp$ is a diffeomorphism of $B_R$ onto 
$D_R$, and let $f\in C_R^\infty (U/K)^K$.
Then $\gL^+_K(U)\ni \mu\mapsto \widehat{f}(\mu)\in \C$ extends to
a holomorphic function on $\fb_\C^*$ of exponential type $R$
such that for all $\lambda \in \fb_\C^*$ and all $w\in W$ we have:
$$\widehat{f}(\lambda )= \widehat{f}(w(\lambda +\rho)-\rho)\, .$$
\end{theorem}
\begin{proof} 
Observe that $\delta(b^{-1})=\delta(b)$ because the 
multiplicities are even.
The formula for $\widehat{f}$ follows then directly from Corollary
 \ref{cor:psi}, formula (\ref{def-Four}) and the $W$-invariance 
of $D$ and $f$. 
For the second part, we note that, by Theorem \ref{th-anaext}, for every fixed
$b \in D_R\cap B$ the map
$$\fb^*_\C\ni \lambda \mapsto \delta (b^{-1})\psi_\lambda (b^{-1})\in \C$$
is holomorphic. Hence, as we are integrating over a compact set,
the map
$$\fb_\C^* \ni \lambda \mapsto 
c\int_{D_R} f(b)\psi_\lambda (b^{-1})\delta (b)\, db\in\C$$
is holomorphic. By Theorem \ref{th-anaext}, we have for the holomorphic
extension:
\begin{equation} \label{eq:extFour} 
d(\lambda ) \widehat{f}(\lambda )=c |W| \; 
\int_{D_R} [D^*f](b) b^{-(\lambda +\rho)}\,db\,.
\end{equation}
The statement on the exponential type $R$ of $\widehat{f}$ follows 
now using Lemma \ref{le-exp} because the right hand side of (\ref{eq:extFour})
is a Fourier transform for the torus $B$.
\end{proof}

\subsection{The local Paley-Wiener Theorem}
In this subsection we prove the non-trivial part
of the local Paley-Wiener Theorem for compact
symmetric spaces with even multiplicities. But first
let us introduce some notation.
Let $\eta: \widetilde{U^\prime}\to U^\prime$ be the universal
covering of $U^\prime$, and set $\widetilde{U}=T\times \widetilde{U^\prime}$.
Let $\exp_{\widetilde{U}}:\fu \to \widetilde{U}$ be
the exponential map, and let
$\eta_1 :=\id\times \eta : \widetilde{U}\to U$.
Then $\exp_U=\eta_1\circ \exp_{\widetilde{U}}\,.$ 
Set $\widetilde{K}:=\exp_{\widetilde{U}}\mathfrak{k}$. Then 
$\eta_1(\widetilde{K})=K$. Since $K \cap T=\{e\}$, we have 
$\widetilde{K}\subset \widetilde{U}'$. The symmetric space 
$\widetilde{U}'/\widetilde{K}$ is the universal covering manifold
of $U'/K$. Let $\widetilde{d}$ denote the Riemannian metric on 
$\widetilde{U}/\widetilde{K}$ induced by the fixed $U$-invariant 
inner product on $\fu$. 
The induced map $\eta_1:\widetilde{U}/\widetilde{K}\to U/K$ is a
local isometry.
Let $\widetilde{x}_0=\{\widetilde{K}\}$ be the base point in 
$\widetilde{U}/\widetilde{K}$. 
Setting $\widetilde{D}_R:=\{\widetilde{x} \in \widetilde{U}/\widetilde{K}
\mid
\widetilde{d}(\widetilde{x},\widetilde{x}_0)\leq R\}$, we have
$\eta_1(\widetilde{D}_R)=D_R$. 
Finally, let $\Exp_{\widetilde{U}}:\fq\to \widetilde{U}/\widetilde{K}$
be defined by 
$\Exp_{\widetilde{U}}(X):=\exp_{\widetilde{U}}X \cdot \widetilde{x}_0$.
Then $\Exp_U = \eta_1 \circ \Exp_{\widetilde{U}}$.

\begin{definition}\label{def:small}
We say that $R>0$ is \emph{small} if the following two 
conditions are satisfied:
\begin{enumerate}
\renewcommand{\labelenumi}{(\theenumi)}
\item $\Exp_{\widetilde{U}}:B_R \to \widetilde{D}_R$ is a diffeomorphism;
\item $\eta_1^{-1}(D_R)$ is a disjoint union of copies 
diffeomorphic to $D_R$ under $\eta_1$.
\end{enumerate}
\end{definition}

Notice that, if $R$ is small, then $\eta_1$ gives a diffeomorphism of
$\widetilde{D}_R$ onto $D_R$. Moreover, in this case, the restriction of 
$\Exp$ is a diffeomorphism of $B_R$ onto $D_R$.

We underline that we are employing the following notion of diffeomorphism
for closed subsets: 
If $C$ and $C'$ are closed subsets of manifolds 
$M$ and $M'$, respectively, then a map $\vartheta:C\to C'$ is said
to be a diffeomorphism if there exists open sets $U$ in $M$ and $U'$ in
$M'$ with $C \subset U$ and $C' \subset U'$, so that $\vartheta:U \to
U'$ is a diffeomorphism. 
According to Theorem \ref{th-spher}, we have
\begin{equation}
\gL_K^+(U)\subseteq
\gL^+_{\widetilde{K}}(\widetilde{U})=\gL^+:=
\left\{\mu\in i\gG^*\oplus \fa^*\mid \forall \alpha\in\Delta^+\, :\,
\mu_\alpha:=\frac{\inner{\mu}{\alpha}}{\inner{\alpha}{\alpha}} 
\in\mathbb{N}_0\right\}\,.
\end{equation}
 We set
\begin{equation}\label{eq-gL}
\gL :=
\{\mu\in i\gG^*_0\oplus \fa^*\mid \forall \alpha \in \gD\, :\,
\mu_\alpha \in\mathbb{\Z}\}\, .
\end{equation}
The decomposition
$$(\mathcal{L},L^2(U/K))=\bigoplus_{\mu \in \gL^+_K(U)}(\pi_\mu,
V_{\pi_\mu})\,,$$
where $\mathcal{L}$
stands for the left-regular representation
of $U$ in $L^2(U/K)$ and $V_{\pi_\mu}$ denotes the Hilbert space of
$\pi_\mu$, implies the following lemma.

\begin{lemma}\label{le-Ktypes}
A function $f\in C(\widetilde{U}/\widetilde{K})$
is of the form $g\circ \eta_1$ for some $g\in C(U/K)$
if and only if $\widehat{f}(\mu )=0$ for all
$\mu\in \gL^+ \setminus
\gL^+_K(U)$. If f has support in $\widetilde{D}_R$, then $g$ has support
in $D_R$.
\end{lemma}

\begin{definition}[Local Paley-Wiener space]  \label{def:PW}
We shall denote by $\mathrm{PW}_{R}(\mathfrak{b}^{\ast})$ the space of
holomorphic functions of exponential type $R$ satisfying 
 $F(w(\lambda+\rho)-\rho
)=F(\lambda)$ for all $\lambda\in\mathfrak{b}^{\ast}_\C$ and
all $w\in W$. Furthermore, we set
\begin{equation}\label{def-PW} 
\mathrm{PW}_R(\fb^*,U):=
\{F\in \mathrm{PW}_R(\fb^*)\mid \forall \mu \in \gL^+\setminus
\gL^+_K(U)\, :\, F(\mu )=0\}\, .
\end{equation}
We call $\mathrm{PW}_R(\fb^*,U)$ the \textit{local Paley-Wiener space} on
$\fb_\C^*$.
\end{definition}

\begin{theorem}[Local Paley-Wiener Theorem]\label{th-pw}
Suppose that $R>0$ is small (according to Definition \ref{def:small}).
Then the Fourier transform $\mathcal{F}$ is a bijection of $C_{R}^{\infty
}(U/K)^K$ onto $\mathrm{PW}_{R}(\mathfrak{b}^{\ast},U)$.
\end{theorem}

We have already seen that the Fourier transform maps
$C_R^\infty (U/K)^K$ injectively into $\mathrm{PW}_R(\fb^*)$, so we
only have to show the surjectivity. Given $F\in \mathrm{PW}_R(\fb^*)$,
then, by the inversion formula in Theorem \ref{th-finv}, we
have to define
\begin{equation}\label{eq-f}
f=\sum_{\mu\in \gL^+_K(U)}d (\mu ) F(\mu )\psi_\mu\,.
\end{equation}
We must show:
\begin{enumerate}
\item $f$ is smooth and $K$-invariant;
\item $\widehat{f}=F$;
\item $\Supp (f)\subseteq D_R$.
\end{enumerate}

We start with some necessary preliminaries. 
Let $X\in \fu$ and $\mu\in \gL^+_K(U)$.
Denote by $\pi_\mu^\infty(X)$ the bounded linear map defined on 
$V_\mu$ by
$$\pi^\infty_\mu (X)v
:=\lim_{t\to 0}\frac{\pi_\mu (\exp tX)v-v}{t}\,,\qquad v\in V_\mu\, .$$
We can extend $\pi_\mu^\infty$ to all of $\fu_\C$ by complex linearity.

\begin{lemma}\label{le-norm} Let $X\in \fu$ and $\mu\in \gL^+_K(U)$. Then
$$\|\pi_\mu^\infty(X)\|\le \|\mu \|\|X\|\, .$$
\end{lemma}
\begin{proof}
Notice that, if $X\in \fu$, then there
exists $k\in U$ such that $\Ad (k)X\in \fc$, where $\fc$
is the Cartan subalgebra from Section 2. Furthermore
$$\|\pi^\infty_\mu (\Ad (k)X)\|=\|\pi_\mu(k)\pi^\infty_\mu
(X)\pi_\mu(k^{-1})\|
=\|\pi_\mu^\infty (X)\|\, .$$
We can therefore assume that $X\in \fc$.

Denote the set of roots of $\fc_\C$ in $\fu_\C$ by
$\gD (\fc)$, the set of positive roots from
Section 2 by $\gD^+(\fc)$ and the corresponding
set of simple roots by $\Sigma (\fc)$. Finally let $W(\fc)$
denote the corresponding Weyl group. 
As $\pi_\mu$ extends to a representation
of $U_\C$ it is easily seen that $\|\pi_\mu^\infty
(w(X))\|=\|\pi_\mu^\infty (X)\|$ for
all $w\in W(\fc)$ and $X\in\fc$. Let
$$i\fc^+=\{X\in i\fc\mid \forall \alpha \in \gD^+(\fc)\, :\, \alpha (X)>0\}\,
.$$
Then, if $X\in i\fc$, there exists $w\in W$ such
that $w(X)\in \overline{i\fc^+}$. Let $w_0$ be
the longest element in $W (\fc)$.
Then there exists a
orthogonal basis $v_\lambda$ 
consisting of 
weight vectors for $\fb$, i.e.\ for all $X\in \fb$ we have
$$\pi^\infty_\mu (X)v_\lambda = \lambda (X)v_\lambda\, .$$
Furthermore, each weight is of the form
$$\lambda = \mu-\sum_{\alpha\in \Sigma (\fc)}n_\alpha \alpha
=w_0\mu +\sum_{\alpha\in \Sigma (\fc)}k_\alpha \alpha$$
for some $n_\alpha,k_\alpha\in\N_0$.
If $X\in i\fc^+$ we get:
$$\mu (X)\ge \mu (X)-\sum n_\alpha \alpha (X)=
(w_0\mu) (X)+\sum k_\alpha \alpha(X)\ge (w_0\mu)(X)\, .$$
As $\|\mu \|=\|w_0\mu \|$ it follows that
$$|\lambda (X)|\le \|\mu \|\|X\|$$
and hence for $X\in \fc$:
$$ |\lambda (X) |=|\lambda (iX)|\le \|\mu\|\, \|X\|\,.$$
The claim in the Lemma therefore follows.
\end{proof}

\begin{lemma}\label{le-pw3}
Let $F\in \mathrm{PW}_R(\fb^*,U)$ and define $f$ by (\ref{eq-f}). Then
$f\in C^\infty (U/K)^K$, and $\widehat{f}(\mu )=F(\mu)$ for
all $\mu\in \gL^+_K(U)$.
\end{lemma}
\begin{proof} Let $\mu \in \gL^+_K(U)\subseteq i\fb^*$.
By  Theorem \ref{th-cfunctdim} we have that
$d(\lambda)$ is a polynomial of degree
$L=2m |\gD^+|$. Furthermore,
$$|\psi_\mu (x)|=|(\pi_\mu (x) w_\mu ,w_\mu)|\le 1\, .$$
It follows that for each $N\in \N$, there exists a constant $C>0$ such that
$$|d (\mu )F(\mu )\psi_\mu (b)|\le C_N(1+\|\mu \|)^{-N}\, .$$
By choosing $N$ large enough, it follows that the series (\ref{eq-f})
converges uniformly. Hence $f$ is continuous. As each $\psi_\mu$ is
$K$-bi-invariant, we deduce that $f$ is $K$-invariant.
Choose a basis $X_1,\ldots ,X_k$ of $\fu$ with $\|X_j\|=1$ for all $j$, and
let
$$g_\mu (t_1,\ldots ,t_k):=\psi_\mu (\exp (t_1X_1)\ldots \exp (t_kX_k))\, .$$
For a fixed $j$, let $g_1=\exp (t_1X_1)\ldots \exp (t_{j-1}X_{j-1})$
and $g_2=\exp (t_{j+1}X_{j+1})\ldots \exp (t_{k}X_{k})$. Then
\begin{eqnarray*}
\bigg|\frac{\partial \, }{\partial t_j}g_\mu  (t_1,\ldots ,t_k)\bigg|
&=&|(\pi_\mu (g_1)\pi_\mu^\infty (X_j)
\pi_\mu (g_2)w_\mu ,w_\mu)|\\
&\le& \|\pi_\mu^\infty (X_j)\|\\
&\le&\|\mu \|
\end{eqnarray*}
by Lemma \ref{le-norm}. Iteration shows that for any multi-index 
$\alpha\in \N_0^k$,
we have
$$ |D^\alpha g_\mu (t_1,\ldots ,t_k)|\le \|\mu \|^{|\alpha|}$$
where $|\alpha|=\alpha_1+\ldots +\alpha_k$.
It follows, as above, that the series $\sum_\mu d(\mu )
F(\mu )D^\alpha g_\mu (\exp(t_1X_k)\ldots
\exp (t_kX_k))$ converges uniformly. Hence $f$ is smooth. 
In particular we have $f\in L^2(U/K)^{K}$, and therefore
$$\sum_{\mu \in \gL^+_K(U)}d(\mu )\widehat{f}(\mu)\psi_\mu
=f=\sum_{\mu \in \gL^+_K(U)}d(\mu )F(\mu) \psi_\mu$$
in $L^2(U/K)$. Taking the inner product with $\psi_\mu$, we see that
$\widehat{f}(\mu )=F(\mu)$ for all $\mu \in \gL^+_K(U)$.
\end{proof}

We will now show that $\Supp (f)\subseteq D_R$. For this, it is
enough to show that $\Supp (\delta f)\subseteq D_R$.
\begin{lemma}\label{le-sumforf} Let $R$ be
small according to Definition \ref{def:small}. 
Let $F\in \mathrm{PW}_R(\fb^*,U)$, and define $f$ by (\ref{eq-f}). 
Then for $b\in B$:
$$\delta (b)f(b)=D\left(\sum_{\mu\in\gL^+_K(U)} 
F(\mu )\sum_{w\in W}b^{w(\mu +\rho)}\right)\, .$$
\end{lemma}
\begin{proof}
If $b\in B$, then by the proof of Lemma \ref{le-pw3}, we have
$$\sum_{\mu\in\gL^+_K(U)} F(\mu )\sum_{w\in W}Db^{w(\mu +\rho)}
= D\left(\sum_{\mu\in\gL^+_K(U)} F(\mu )\sum_{w\in W}b^{w(\mu +\rho)}\right)\,
.$$
Hence, for all $b\in B$ with $\delta (b)\not= 0$, we get
\begin{eqnarray*}
\delta (b) f(b)&=&\delta (b) \sum_{\mu\in\gL^+_K(U)}d(\mu ) F(\mu )
\psi_\mu (b)\\
&=&\sum_{\mu\in\gL^+_K(U)} F(\mu )\sum_{w\in W}Db^{w(\mu +\rho)}\\
&=& D\left(\sum_{\mu\in\gL^+_K(U)} F(\mu )\sum_{w\in W}b^{w(\mu
    +\rho)}\right)\, .
\end{eqnarray*}
As both sides are continuous in $b$, it follows that this holds on
$\overline{\{b\in B\mid \delta (b)\not= 0\}}=B$.
\end{proof}

Before finishing the proof of the main theorem, we need 
the following well-known lemma.
Recall that, for $\mu \in i\fb^*$, we have introduced 
the notation $\chi_\mu (b):=b^\mu$, provided $b^\mu$
is defined for all $b\in B$.

\begin{lemma}\label{le-GammaUK} Let
$\Gamma_1:=\{X\in\mathfrak{b}_1\mid\exp_{\widetilde{U^\prime}}(X)\in
\widetilde{K}\}$.
Then
$$\gG_1=\{X\in \fb_1\mid \forall \mu \in \gL \, :\, \mu (X)\in 2\pi i\Z\}\,
.$$
Furthermore, if $\Gamma =\Gamma_0\oplus \Gamma_1$, then
$i\Gamma^*=\gL$ and
the map
$$\gL \ni \mu\mapsto \chi_\mu\in \widehat{\fb/\gG} $$
is a bijection.
\end{lemma}

\begin{proof} See the proof of Lemma 4.1, p.\ 535, in \cite{He84}.
\end{proof}

\begin{lemma}\label{le-sums}
Let $F\in\mathrm{PW}_{R}(\mathfrak{b}^{\ast})$ and $b\in B$. Then
\begin{equation*}
D\left(\sum_{\mu\in\Lambda^{+}}
F(\mu)\sum_{w\in W}b^{w(\mu+\rho)}\right)=
D\left(\sum_{\mu\in\Lambda}
F(\mu-\rho)b^{\mu}\right)\,.
\end{equation*}
\end{lemma}

\begin{proof}
Set $G(\mu):=F(\mu-\rho)$. Then $G$ is of exponential type $R$. 
It follows, as in the proof  of Lemma \ref{le-pw3},
that $\sum_{\mu\in\Lambda}G(\mu)b^{\mu}$   defines  a
smooth function on $B$. From $F(w(\lambda +\rho)-\rho)=F(\lambda)$ 
we obtain that $G(w(\lambda+\rho))=G(\lambda+\rho)$ for all 
$\lambda \in \mathfrak{b}_{\mathbb{C}}^{\ast}$ and $w\in W$.
Finally, part (2) of Lemma \ref{le-dim} implies that
$\mu\mapsto\mu+\rho$ is a bijection on $\Lambda$.
As $\Lambda^+$ is a fundamental domain for the action of $W$ on 
$\Lambda$, we have:
\begin{eqnarray*}
\sum_{\mu\in\Lambda }G(\mu)b^{\mu} & =&\sum_{\mu\in\Lambda }G(\mu
+\rho)b^{\mu+\rho} \\
& =&\sum_{\mu+\rho\in\Lambda^{+}}F(\mu)\frac{1}{|W^{\mu+\rho}|}\sum_{w\in
W}b^{w(\mu+\rho)}
\end{eqnarray*}
where $W^{\mu+\rho}:=\{w\in W\mid w(\mu+\rho)=\mu+\rho\}$.
 
For the final step, assume first that
$\mu+\rho\in\Lambda^{+}$, but $\mu\not \in \Lambda^{+}$. Then there
is a simple root $\beta\in\Delta^{+}$ such that $\inner{\mu}{\beta}<0$. 
As $\mu\in \Lambda$, it follows that $\mu_\beta\in \Z$.
In particular, $\mu_\beta \le -1$. Since
$\inner{\mu+\rho}{\beta}\ge 0$ we have, 
using part (1) of Lemma \ref{le-dim},
\begin{equation*}
0\leq \frac{\inner{\mu+\rho}{\beta}}{\inner{\beta}{\beta}}
=(\mu+\rho)_\beta \leq -1+m \,.
\end{equation*}
By Corollary \ref{co-zero} it follows that
\begin{equation*}
D\left(\sum_{w\in W}b^{w(\mu+\rho)}\right)=0\,.
\end{equation*}
Finally, if $\mu+\rho\in\Lambda_{K}^{+}$ and $\mu\in\Lambda_{K}^{+}$, 
then $W^{\mu+\rho}=\{e\}$. The claim thus follows.
\end{proof}

\begin{lemma}
\label{le-support} 
Suppose that $R>0$ is small in the sense of 
Definition \ref{def:small}. 
Let $F\in\mathrm{PW}_{R}(\mathfrak{b}^{\ast})$.  
Define $h:\mathfrak{b}\rightarrow\mathbb{C}$ by
\begin{equation*}
h(X)=\sum_{\mu\in\Lambda}F(\mu-\rho)e^{\mu(X)}~.
\end{equation*}
Then $h$ is a $\Gamma$-periodic smooth function on $\fb$ and 
$\mathrm{Supp}(h)\subseteq B_R+\gG$.
\end{lemma}

\begin{proof}
It follows from Lemma \ref{le-GammaUK} and from the proof of Lemma 
\ref{le-pw3} that $h$ is smooth and $\Gamma $-periodic. 
It therefore defines a smooth function on the
abelian group $\mathfrak{b}/\Gamma$.
Hence
\begin{equation*}
F(\mu-\rho)=\mathrm{vol}(\mathfrak{b}/\Gamma)^{-1}
\int_{\mathfrak{b}/\Gamma}h(X)e^{-\mu(X)}~dX=\widehat{h}(\mu)
\end{equation*}
By the classical Paley-Wiener Theorem there is a
$g\in C_{R}^{\infty}(\mathfrak{b})$
such that $\widehat{g}(\lambda)=F(\lambda-\rho)$, 
Here the Fourier
transform of $g$ is defined by
\begin{equation*}
\widehat{g}(\lambda)=\frac{1}{(2\pi)^{n}}
\int_{\mathfrak{b}}g(X)e^{-\lambda(X)}~dX~,\quad\lambda\in
i\mathfrak{b}^{\ast}
\end{equation*}
where $n=\dim \fb$, as before.
We claim that there exists
a constant $\gamma \not= 0$ such that
$$\sum_{Y\in\Gamma}g(X+Y)=\gamma\, h(X)\, .$$
Indeed, let
$$G(X)=\sum_{Y\in\Gamma}g(X+Y)\, .$$
Then $G$ is $\Gamma$-periodic and%
\begin{eqnarray*}
\widehat{G}(\mu) & =&\mathrm{vol}(\mathfrak{b}/\Gamma)^{-1}
\int_{\mathfrak{b}/\Gamma}G(X)e^{-\mu
(X)}~dX \\
& =&\mathrm{vol}(\mathfrak{b}/\Gamma)^{-1}
\int_{\mathfrak{b}}g(X)e^{-\mu(X)}~dX \\
& =&\mathrm{vol}(\mathfrak{b}/\Gamma)^{-1}
(2\pi)^{n}\widehat{g}(\mu) \\
& =&\mathrm{vol}(\mathfrak{b}/\Gamma)^{-1}
(2\pi)^{n}F(\mu-\rho)~ \\
& =&\mathrm{vol}(\mathfrak{b}/\Gamma)^{-1}
(2\pi)^{n}\widehat{h}(\mu)~.
\end{eqnarray*}
But this implies that
$$G =\mathrm{vol}(\mathfrak{b}/\Gamma)^{-1}(2\pi)^{n}h\, .$$
Observe that $(B_R+\gamma_1)\cap (B_R+\gamma_2)=\emptyset$
if $\gamma_1,\gamma_2\in \gG$ and $\gamma_1\not=\gamma_2$. Hence
$\Supp (g)\subseteq B_R$ implies 
that $\mathrm{Supp}(G)\subseteq B_R+\gG$. 
The same must therefore hold for $h$.
\end{proof}

We now finish the proof of the local Paley-Wiener Theorem by proving
the following lemma:
\begin{lemma}
Assume that $R>0$ is small in the sense of Definition \ref{def:small}
and that $F\in\mathrm{PW}_{R}(\mathfrak{b}^{\ast},U)$. Then 
there exists a $f\in C_{R}^{\infty}(U/K)^{K}$ such that 
$\widehat{f}(\mu)=F(\mu)$ 
for all $\mu\in\Lambda_{K}^{+}$. Hence the spherical Fourier transform
$\cF : C_R^\infty (U/K)^K\to \mathrm{PW}_R(\fb^\ast,U)$ is
surjective.
\end{lemma}

\begin{proof} By Lemma \ref{le-Ktypes} we can assume that
 $U=\widetilde{U}$ and $K=\widetilde{K}$. Hence $\Lambda^+_K=\Lambda^+$.
Define a smooth $K$-invariant function $f$ on $U/K$ by%
\begin{equation*}
f(x)=\sum_{\mu\in\Lambda^{+}}d(\mu)F(\mu)\psi_{\,\mu}(x) \, .
\end{equation*}
Then $\widehat{f}(\mu )=F(\mu )$ for all $\mu\in \gL^+$,
cf.\ Lemma \ref{le-pw3}.
As already observed, it sufficed to prove that $\Supp(f|_B) \subseteq D_R$,
which is equivalent to the condition $\Supp(\delta f|_B) \subseteq D_R$.
By Lemma \ref{le-sumforf} and Lemma \ref{le-sums} we have for $X \in \fb$:
\begin{equation}\label{eq:deltafDh}
(\delta f)(\exp X)= D\Big( \sum_{\mu\in\Lambda^+}F(\mu)
\sum_{w\in W}e^{w(\mu+\rho)(X)}\Big)
= D\Big(\sum_{\mu\in \gL}F(\mu -\rho)e^{\mu(X)} \Big)
= D h(X)
\end{equation}
where $h$ is the function defined in Lemma \ref{le-support}.
Here $Dh$ is defined locally by $Dh:=D(h \circ \exp^{-1})$.
According to Lemma \ref{le-support}, we know that
$h$ (and hence $Dh$) has support in $B_R+\Gamma$. Thus $\delta\, f$
has support in $D_R$.
\end{proof}

We conclude this section by proving two integral formulas for the
smooth functions on $U/K$ with ``small support''. The first formula,
which can be deduced from the proof of the local Paley-Wiener theorem, 
will play a decisive role in proving the validity of Huygens' 
principle on $U/K$.

\begin{corollary} \label{cor:integralformula}
Suppose that $R>0$ is small in the sense of Definition \ref{def:small}
and that $f\in C_{R}^{\infty}(U/K)^{K}$.
Then the following integral formulas hold on $B$: 
\begin{equation}\label{eq:intform}
\delta(b)f(b)=D\left( \int_{i\fb^*} \widehat{f}(\lambda-\rho) b^\lambda\; 
d\lambda\right)\,,
\end{equation}
where $D$ is the differential operator of Theorem 
\ref{th-anaext}, and
\begin{equation*}
f(b)=\frac{1}{|W|} \; \int_{i\fb^*} \widehat{f}(\lambda-\rho) 
\varphi_\lambda(b) \,d(\lambda) \, d\lambda\,.
\end{equation*}
\end{corollary}
\begin{proof}
Suppose $b=\exp X \in D_R \cap B$. Then by (\ref{eq:deltafDh}), we have
$\delta(b)f(b)=Dh(X)$. The function $h$ is determined by the proof of
Lemma \ref{le-support} with $F=\widehat{f}$. Keeping the notation of that
proof, we obtain $h(X)=g(X)$ because $X \in B_R$, and 
$\widehat{g}(\l)=\widehat{f}(\l-\rho)$ for all $\l \in \fb_\C^*$. Hence
$$
h(X)=g(X)=\int_{i\fb^*} \widehat{g}(\l) e^{-\l(X)}\;d\l =
\int_{i\fb^*} \widehat{f}(\l-\rho) e^{-\l(X)}\;d\l=
\int_{i\fb^*} \widehat{f}(\l-\rho) b^{-\l}\;d\l\,.
$$
This proves (\ref{eq:intform}) because both $f|_B$ and 
$\int_{i\fb^*} \widehat{f}(\l-\rho) b^{-\l}\;d\l$ are supported in 
$\exp B_R$. The last formula follows then immediately
from Theorem \ref{th-anaext} and the $W$-invariance of $\l \mapsto 
\widehat{f}(\l-\rho)$.
\end{proof}

\section{The local Huygens' principle for compact symmetric spaces with even
multiplicities}\label{section:Huygens}

\noindent
Let $L_X$ denote the Laplace-Beltrami operator on 
a Riemannian symmetric space $X$ of the noncompact or compact type.
The \emph{modified wave equation} 
on $X$ is the partial differential equation
\begin{equation}\label{eq:wave}
(L_X\pm \|\rho\|^2)u=u_{tt}\,,
\end{equation}
where $u=u(x,t)$ is a function of $(x,t) \in X \times I$ and $I\subseteq \R$
is an interval containing $0$. The sign in front of $\|\rho\|^2:=
\inner{\rho}{\rho}$ has to be chosen $+$ if $X$ is of the noncompact type,
and $-$ if $X$ is of the compact type.
 
Let $f \in C^\infty_c(X)$ be fixed. Huygens' principle concerns specific 
support properties for the smooth solution $u$ of (\ref{eq:wave}) which 
satisfies the Cauchy conditions
\begin{equation}\label{eq:cauchy}
\begin{array}{rl}
u(x,0)&\!\!\!=0,  \\   
u_t(x,0)&\!\!\!=f(x).
\end{array}
\end{equation}
Recall that solving a Cauchy problem with initial conditions
$u(x,0)=g(x)$, $u_t(x,0)=f(x)$, where $f, g \in C^\infty_c(X)$ are arbitrary,
can always be reduced to solving a Cauchy problem with initial conditions
of the form (\ref{eq:cauchy}), and that support properties like Huygens'
principle reduce at the same time.  
Indeed, if $u_i$ for $i=1,2$
is the solution to $Lu_i=(u_i)_{tt}$ (for $L$ the operator on the 
left in (\ref{eq:wave}))
with Cauchy data $(0,f_i)$, then $u:=u_2+(u_1)_t$ is
a solution with Cauchy data
$$
\begin{array}{rl}
u(x,0)&\!\!\!=f_1(x), \\
u_t(x,0)&\!\!\!=f_2(x)+(u_1)_{tt}(x,0)=f_2(x)+(Lu_1)(x,0)=f_2(x).
\end{array}
$$

Moreover, if $u(x,t)$ is the solution 
corresponding to (\ref{eq:cauchy}), then $u(x,-t)$ corresponds to the initial
conditions $u(x,0)=0$, $u_t(x,0)=-f(x)$. This allows us to restrict our 
analysis to values $t\ge 0$. Finally, the general case can be reduced to 
the $K$-invariant one.
We shall therefore assume in the following that $f \in C^\infty_c(X)^K$. 
In this case, the solution $u$ will be a $K$-invariant function 
of the variable $x \in X$.

The property that the support of the solution $u$ is compact is stated by the 
\emph{principle of finite propagation speed}. This principle holds, 
more generally, for solutions of wave equations on arbitrary 
Riemannian manifolds.
(See e.g.\ \cite{Fr}, Ch.\ 5.) In our setting, it corresponds
to the following lemma. We recall the notation $D_R:=\{x\in X\mid
d(x,x_0)\leq R\}$ for the closed ball of center
$x_0$ and radius $R$.  

\begin{lemma}[Finite propagation speed]
\label{lemma:speed}
Let $u \in C^\infty(X\times I)$ be solution
to the Cauchy problem \straightparen{\ref{eq:cauchy}}. Let $\e>0$ and let 
$t\in I \cap (0,+\infty)$.
Suppose that the Cauchy datum $f$ in \straightparen{\ref{eq:cauchy}}
satisfies $\Supp(f)\subseteq D_\e$. Then
$\Supp\big(u(\cdot,t)\big) \subset D_{\e+t}$.
\end{lemma}
 
For positive values of $t \in I$, the solution $u(x,t)$ to the Cauchy 
problem is therefore supported inside the positive cone
\begin{equation}\label{eq:econe} 
C_\e:=\{(x,t) \in X \times [0,\infty)\mid d(x,x_0)\leq \e+t\}.
\end{equation}

Let $\e$, $t$, $f$ and $u$ be as in Lemma \ref{lemma:speed}.
We say that the \emph{\straightparen{strong} Huygens' principle} holds
provided 
$\Supp\big(u(\cdot,t)\big) \subseteq  \{x\in X\mid t-\e\leq d(x,x_0)\leq 
t+\e\}$. 

Thus, when the strong Huygens' principle holds, the solution of 
$u(x,t)$ to the Cauchy problem is supported for positive $t \in I$ 
inside the conical shell
\begin{equation} \label{eq:eshell}
S_\e:=\{(x,t) \in X \times [0,\infty)\mid t-\e\leq d(x,x_0)\leq t+\e\}.
\end{equation}

Huygens' principle holds true (at least for small values of $t$) 
for the wave equation on odd dimensional Riemannian symmetric spaces
$X$ of the noncompact or compact type for which all root multiplicities are
even. (See the references in the introduction.) As can be verified from 
the tables reported in Subsection \ref{subsection:classification},
$\dim(U/K)$ is odd if and only if $\dim(\fb)=\rank(U/K)$ is odd.

The main result of this section is the following local version of Huygens'
principle on symmetric spaces of the compact type.

\begin{theorem}\label{thm:HP}
Let $U/K$ be a Riemannian symmetric space with all even multiplicities.
Let $R>0$ be small according to Definition \ref{def:small}. Let
$0<\e<R$, and let $f \in C^\infty_\e(U/K)^K$. 
Assume that $U/K$ is odd dimensional \straightparen{that is, 
$\rank(U/K)=\dim\fb$ is odd}. 

Suppose that $u(x,t)$ is a smooth solution of Cauchy's problem
\begin{align}\label{eq:cauchyUK}
(L_{U/K}-\|\rho\|^2)u&=u_{tt}\,, \notag\\
u(x,0)&=0,  \\   
u_t(x,0)&=f(x)\,. \notag
\end{align}
Then the following properties are satisfied:
\begin{enumerate}
\renewcommand{\theenumi}{\alph{enumi}}
\renewcommand{\labelenumi}{(\theenumi)}
\item \emph{(Local exponential Huygens' principle)} 
There is a constant $C>0$ so that
for all $(x,t) \in U/K\times [0,R-\e]$ 
and all $\gamma \in [0,\infty)$ we have
\begin{equation} \label{eq:expHP}
|\delta(x)u(x,t)|\leq C e^{-\gamma(t-d(x,x_0)-\e)}\,.
\end{equation}
Here $\delta$ denotes the $K$-invariant extension to $U/K$ of the 
$W$-invariant function $\delta$ defined in \straightparen{\ref{eq:delta}}.
\item \emph{(Local strong Huygens' principle)}
$$\Supp(u) \cap ( U/K \times [0,R-\e])= 
\Supp(u) \cap ( D_R \times [0,R-\e]) \subseteq S_\e\,,$$
where $S_\e$ denotes the $\e$-shell \straightparen{\ref{eq:eshell}}.

\begin{center}
\includegraphics{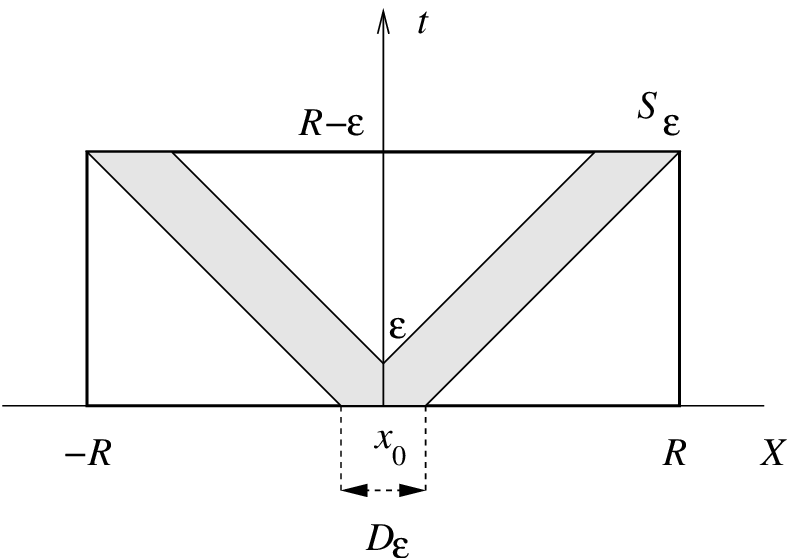}
\end{center}
\item 
Suppose $\dim(U/K) \geq 3$. Let $D$ be the differential operator of 
Theorem \ref{th-anaext}.
Then for all $b=\exp X \in B$ and $t \in [0,R-\e]$ 
the smooth solution $u(b,t)$ to \straightparen{\ref{eq:euclCauchy}} is given
by
formula 
\begin{equation}
\label{eq:uodd}
\delta(b)u(b,t)=\frac{\Omega_n/2}{[(n-3)/2]! \Omega_{n-1}}\; 
D\Big( \frac{\partial}{\partial(t^2)}\Big)^{(n-3)/2} 
  \big( t^{n-2} (M^tf)(X)\big)\,.
\end{equation}
Here 
\begin{equation}\label{eq:M}
(M^rf)(X):=\frac{1}{\Omega_{n-1}(r)} \int_{S_r(X)} f(s)\,d\sigma(s)\,
\end{equation}
is the mean value of a function $f:\fb\to \C$ on the Euclidean sphere 
  $S_r(X):=\{X\in \fb\mid 
\|X\|=r\}$ in $\fb \cong \R^n$ with respect to the 
${\rm O}(n)$-invariant surface measure $d\sigma$. Moreover, 
$\Omega_{n-1}(r)$ denotes the surface area of $S_r(x)$, and 
$\Omega_{n-1}:=\Omega_{n-1}(1)$.
\end{enumerate}
\end{theorem}

Parts (a) and (c) of  Theorem \ref{thm:HP} seem to be new in the context
of symmetric spaces of compact type. As we shall see in the following, they
both imply the local strong Huygens' principle of Part (b). Another
independent proof of Theorem \ref{thm:HP} will be given in Corollary 
\ref{cor:firstHP}. 

The remainder of this section is devoted to the proofs of the three
parts of 
Theorem \ref{thm:HP}. To underline the various necessary steps, we
have subdivided them into different lemmas and corollaries. 
Before entering the details of the proofs, we remark that, since the solution 
$u(x,t)$ is smooth and $K$-invariant in the $x$-variable, it suffices 
to examine its restriction to $B \times [0,R-\e]$. This will be common to all
three methods which we are going to describe.

Recall that for all $\mu \in \gL_K(U)^+$ we have
\begin{equation}\label{eq:eigenL}
L_{U/K}\psi_\mu=-\inner{\mu+2\rho}{\mu} \psi_\mu\,.
\end{equation}
By Lemma \ref{lemma:speed}, for fixed $t>0$ the solution 
$u(\cdot,t)$ to (\ref{eq:cauchyUK}) is supported inside $D_{t+\e}$.
This allows us to interchange integration 
and differentiation with 
respect the variable $x \in U/K$. Hence, taking the spherical Fourier
transform of (\ref{eq:cauchyUK}) for fixed $t$, we obtain:
\begin{align}\label{eq:cauchyF}
-\|\mu+\rho\|^2\; \widehat{u}(\mu,t)&=\widehat{u}_{tt}(\mu,t), \notag\\
\widehat{u}(\mu,0)&=0,  \\   
\widehat{u}_t(\mu,0)&=\widehat{f}(\mu)\,. \notag
\end{align}

\begin{lemma}
Let $U/K$ be a Riemannian symmetric space with all even multiplicities. 
Let $R>0$ be small according to Definition \ref{def:small}, and let
$0<\e<R$. 
Let $u(x,t)$ be a smooth solution of Cauchy's problem (\ref{eq:cauchyUK})
with Cauchy datum $f \in C^\infty_\e(U/K)^K$.
Then 
\begin{equation}\label{eq:omega}
\widehat{u}(\l-\rho,t)=\widehat{f}(\l-\rho) \; 
\frac{\sin(\inner{\l}{\l}t)}{\inner{\l}{\l}}\,.
\end{equation}
Consequently, for all $(b,t) \in B \times [0,R-\e]$ we have 
 \begin{equation} \label{eq:deltau}
\delta(b)u(b,t)=D\left( \int_{i\fb^*} \widehat{f}(\lambda-\rho)
\frac{\sin(\inner{\l}{\l}t)}{\inner{\l}{\l}}
 b^\lambda\; d\lambda\right)\,.
\end{equation}
\end{lemma}
\begin{proof}
Suppose that $t \in (0,R-\e)$. Then $t+\e<R$ is small, and 
the local Paley-Wiener Theorem \ref{th-pw}
ensures that $\mu \mapsto \widehat{u}(\mu,t)$ extends uniquely to 
$\l\mapsto \widehat{u}(\l,t) \in {\rm PW}_{t+\e}(\fb^*,U)$. Likewise, 
$\widehat{u}_{tt}(\mu,t)$ and  $-\inner{\mu-\rho}{\mu-\rho}
\widehat{u}(\mu,t)$ admit unique holomorphic extensions in 
${\rm PW}_{t+\e}(\fb^*,U)$, respectively to 
$\widehat{u}_{tt}(\lambda,t)$ and $-\inner{\lambda-\rho}{\lambda-\rho}
\widehat{u}(\lambda,t)$. Finally,  $\widehat{u}(\mu,0)$, 
$\widehat{u}_{t}(\mu,0)$ and
$\widehat{f}(\mu)$ extend uniquely to ${\rm PW}_{\e}(\fb^*,U)$. 
By uniqueness, we
conclude that the equations in (\ref{eq:cauchyF}) hold for the holomorphic
extensions. Setting $\omega(\l,t):= \widehat{u}(\l-\rho,t)$, we 
are therefore reduced to the Cauchy problem
\begin{align*}
\omega_{tt}(\l,t)&=-\inner{\l}{\l}\omega(\l,t)\notag\\
\omega(\l,0)&=0  \\   
\omega_t(\l,0)&=\widehat{f}(\l-\rho)\,, \notag
\end{align*}
from which (\ref{eq:omega}) follows.

By Lemma \ref{lemma:speed}, $u(\cdot,t) \in C^\infty_{t+\e}(U/K)^K$. 
Since $t+\e$ is small
according to Definition \ref{def:small}, formula (\ref{eq:deltau})
is then a consequence of (\ref{eq:intform}) and
(\ref{eq:omega}).
\end{proof}

Define $v:\fb \times (-R+\e,R-\e) \to \C$ by
\begin{equation}\label{eq:v}
v(X,t):=\int_{i\fb^*} \widehat{f}(\lambda-\rho)
\frac{\sin(\inner{\l}{\l}t)}{\inner{\l}{\l}}
 e^{\lambda(X)}\; d\lambda\,.
\end{equation}
Then $v(X,t)$ is the solution of the Cauchy problem for the 
wave equation on $\fb \cong \R^n$:
\begin{align}
L_\fb v(X,t)&=v_{tt}(X,t)  \notag\\
v(X,0)&=0 \label{eq:euclCauchy} \\ 
v_t(X,0)&=g(X)\,, \notag
\end{align}
where $g \in C^\infty_\e(\fb)^W$ is the inverse Euclidean Fourier transform 
of $\widehat{f}(\l-\rho)$.

Since $\exp:B_R \to D_R$ is a diffeomorphism and since 
the operator $D$ preserves supports, we have proven the
following corollary, yielding the first proof of the local 
strong Huygens principle of
Theorem \ref{thm:HP}.

\begin{corollary} \label{cor:firstHP}
Suppose $U/K$ is a symmetric space of the compact type with even 
multiplicities. Let $R$ be small according to Definition \ref{def:small},
and let $0<\e<R$.
Then the strong Huygens' principle holds for (\ref{eq:cauchyUK}) 
on $U/K \times [0,R-\e]$ provided it holds for (\ref{eq:euclCauchy}) 
on $\fb \times[0,R-\e]$.
Hence the local strong Huygens' principle holds if $\dim(U/K)$ is odd
\straightparen{i.e. if $\rank(U/K)=\dim \fb$ is odd}.
\end{corollary}

To prove the local exponential Huygens' principle, we apply the procedure
of \cite{BOS95} to the integral appearing at the right-hand side of 
(\ref{eq:deltau}). Our computations are nonetheless easier than those
in that article. Since we only consider the even multiplicity situation,
we can employ our differential operator $D$. This allows us 
to work in a Euclidean setting by replacing the spherical functions appearing
in the integral formulas studied in \cite{BOS95} with exponential functions. 

Let $S$ denote the unit sphere in $i\fb^*$, and, as before, let $n=\dim
\fb=\rank (U/K)$.
With respect to polar coordinates $(\omega,p) \in S \times [0,+\infty)$ 
in $i\fb^*$, we have $d\l=p^{n-1} d\omega dp$.

Setting 
\begin{equation}\label{eq:Psi}
\Psi_\e(p,X):=p^{n-1} \int_S \widehat{f}(p\omega-\rho) e^{p\omega(X)}\,
 d\omega\,,
\end{equation}
we obtain for $b=\exp X \in B$
\renewcommand{\arraystretch}{2.3}
\begin{equation}\label{eq:duPsi}
\begin{array}{rl}
\delta(b) u(b,t)&\!\!\!= D\left(\displaystyle
\int_{i\fb^*} \widehat{f}(\lambda-\rho) \,
\frac{\sin(\inner{\l}{\l}t)}{\inner{\l}{\l}}\; 
 e^{\lambda(X)}\; d\lambda\right) \\
&\!\!\!= D\left( \displaystyle
\int_0^\infty \frac{\Psi_\e(p,X)}{p}\, \sin(pt) \, dp\right)\,.
\end{array}
\end{equation}
\renewcommand{\arraystretch}{1.3}

\begin{lemma}  \label{lemma:Psi}
Suppose $n:=\dim\fb$ is odd. Then the following properties hold. 
\begin{enumerate}
\renewcommand{\theenumi}{\alph{enumi}}
\renewcommand{\labelenumi}{(\theenumi)}
\item 
The function $\Psi_\e(p,X)$ is holomorphic on $\C \times \fb_\C$. It is 
even in $p\in \C$ and $W$-invariant in $X \in \fb$.
Moreover, for every $N \in \N$ there is a constant $K_N>0$ such that
\begin{equation} \label{eq:estPsi}
\big|\Psi_\e(p,X)\big| \leq K_N |p|^{n-1} (1+|p|)^{-N} 
e^{|\im p|(\e+\|X\|)}
\end{equation}
for all $p \in \C$ and $X \in \fb$.
\item
Suppose furthermore that $n\neq 1$. Let $D$ be the differential operator
  of Theorem \ref{th-anaext}.
Then the function $p^{-1}D\Psi_\e(p,X)$
is holomorphic on 
$\C \times \fb_\C$. It is odd in $p\in \C$ and $W$-invariant in $X \in \fb$.
Moreover, for every $N \in \N$ and every compact $Q \subset \fb$  
there is a constant $K_{N,Q}>0$
such that
\begin{equation} \label{eq:estDPsi}
\bigg| \frac{D\Psi_\e(p,X)}{p}\bigg| \leq K_{N,Q} |p|^{n-2} (1+|p|)^{-N} 
e^{|\im p|(\e+\|X\|)}
\end{equation}
for all $p \in \C$ and $X \in Q$.
\end{enumerate}
\end{lemma}
\begin{proof}
As the integrand in (\ref{eq:Psi}) is holomorphic in 
$(p,X) \in \C \times \fb_\C$ and continuous in $\omega \in S$, 
it follows from Morera's theorem that $\Psi_\e$ is holomorphic in
$\C \times \fb_\C$. The fact that $\Psi_\e$ is even in $p$ and 
$W$-invariant in $X$ is a consequence of the ${\rm O}(\fb^*)$-invariance of 
$d\omega$, the $W$-invariance of $\l \mapsto \widehat{f}(\l-\rho)$,
and the fact that $n$ is odd.

Recall the notation (\ref{eq:ReIm}) for the real and imaginary parts
in $\fb_\C^*$. 
If $(p,\omega)\in \C \times S$, then 
$$\|\re(p\omega-\rho)\|= \|\re(p\omega)\|=|\im p|\|\omega\|=|\im p|\,$$
and
$$\re(p\omega(X))=-\im p \;\im \omega(X) \leq |\im p| \|X\|\,.$$
Since $\widehat{f} \in {\rm PW}_\e(\fb^*)$, we therefore obtain for all
$p\in \C$, $\omega \in S$ and $X \in \fb$:
\begin{align*}
\big|\widehat{f}(p\omega-\rho) e^{p\omega(X)}\big|
&\leq C_N (1+\|p\omega-\rho\|)^{-N} e^{\e\|\re(p\omega-\rho)\|}
   e^{\re(p\omega(X))}\\
&\leq C_N (1+|p|)^{-N} e^{\e|\im p|} 
   e^{|\im p|\|X\|}\,,
\end{align*}
from which the estimate (\ref{eq:estPsi}) immediately follows.
Formula (\ref{eq:Psi}) shows then that $p^{-1}\Psi_\e(p,X)$ remains 
holomorphic provided $n>1$. It is odd in $p \in \C$ and $W$-invariant
in $X \in \fb$. The same property holds therefore also for 
$p^{-1}D\Psi_\e(p,X)$ because $D$ is $W$-invariant and has holomorphic
coefficients.
Suppose $n\geq 3$. Differentiation under integral sign gives 
\begin{equation*}
\frac{D\Psi_\e(p,X)}{p}=p^{n-2} \int_S 
\widehat{f}(p\omega-\rho) De^{p\omega(X)}\, d\omega\,.
\end{equation*}
Notice that 
$De^{p\omega(X)}=f(X,\omega,p) e^{p\omega(X)}$
where $f(X,\omega,p)$ is holomorphic in $X \in \fb^*_\C$, polynomial in 
$\omega \in S$ and polynomial in $p \in \C$.
For every compact subset $Q$ of $\fb$ there is a constant $C_K$
so that 
 $$|f(X,\omega,p)|\leq C_K (1+|p|)^s\,$$
where $s=\deg D$ is the polynomial degree of $f$ in the variable $p$.
This, with the same argument used for (\ref{eq:estPsi}), proves the 
estimate (\ref{eq:estDPsi}).
\end{proof}

We now use Lemma \ref{lemma:Psi} to prove exponential estimates for the
solution $u(x,t)$.
Since $\Psi_\e(p,X)/p$ is odd, we obtain from (\ref{eq:duPsi})
that for all $b=\exp X \in B$, $t \in [0,R-\e]$ and $\gamma >0$ 
we have
\renewcommand{\arraystretch}{2.3}
\begin{equation}\label{eq:shift}
\begin{array}{rl}
\delta(b) u(b,t)
&\!\!\!=\dfrac{1}{2}\; 
D\left( \displaystyle\int_{-\infty}^\infty 
\frac{\Psi_\e(p,X)}{p}\; e^{ipt} \, dp\right) \\
&\!\!\!=\dfrac{1}{2}\; 
\displaystyle\int_{-\infty}^\infty 
\frac{D\Psi_\e(p,X)}{p}\; e^{ipt} \, dp \\
&\!\!\!=\dfrac{1}{2}\; 
\left(\displaystyle\int_{-\infty}^\infty 
\frac{D\Psi_\e(p+i\gamma,X)}{p+i\gamma}\; e^{ipt} \, dp\right)e^{-\gamma t}. 
\end{array}
\end{equation}
In the above computations, the differentiation under integral sign
and the shift in the path of integration are justified by the 
estimates of Lemma \ref{lemma:Psi}.
\renewcommand{\arraystretch}{1.3}

\begin{lemma} \label{lemma:expbigrank}
Under the assumptions of Theorem \ref{thm:HP}, 
the local exponential Huygens' principle of Theorem \ref{thm:HP}(b)
holds when $\dim U/K >1$.
\end{lemma}
\begin{proof}
Equation (\ref{eq:shift}) together with estimate (\ref{eq:estDPsi})
give for all $b=\exp X \in B$ and $\gamma \in [0,\infty)$
\begin{align*}
|\delta(b) u(b,t)|&\leq
\frac{1}{2}\; 
\left(\dint_{-\infty}^\infty 
\bigg|\frac{D\Psi_\e(p+i\gamma,X)}{p+i\gamma}\bigg| |e^{ipt}| \, 
dp\right)e^{-\gamma t}\\
&\leq C_N \Big( \int_{-\infty}^\infty (1+|p|)^{-N} \, dp\Big) 
e^{-\gamma(t-\|X\|-\e)}
\end{align*}
Since $\|X\|=d(x,x_0)$ when $x\in K\exp(X)K$, we conclude that
for all $(b,t) \in B \times [0,R-\e]$ and all $t \in [0,\infty)$, we have
$$|\delta(b) u(b,t)|\leq C e^{-\gamma(t-d(x,x_0)-\e)}\,,$$
where $C$ is a positive constant. The inequality then extends by 
$K$-invariance to $U/K$.
\end{proof}

As in \cite{BOS95}, the exponential estimates  
need to be be worked out directly  when $\rank(U/K)=1$. 
In this case the definition of 
$\Psi_\e(p,X)$ simplifies since $S=\{\pm i\}$. 
We shall identify $\fb_\C^*$ 
with $\C$ by $\l \equiv \usfrac{\inner{\alpha}{\l}}{\inner{\alpha}{\alpha}}$. 
As $\l\mapsto \widehat{f}(\l-\rho)$ is even, it follows that 
\begin{align*}
\Psi_\e(p,X)&=\widehat{f}(ip-\rho)e^{ipX} +\widehat{f}(-ip-\rho)e^{-ipX}\\
           &=2\widehat{f}(ip-\rho)\cos(pX)\,.
\end{align*}
In the rank one case we can write the operator $D$ in the form 
$D=D'(\usfrac{d}{dX})$, where the $D'$ is an odd differential operator with 
holomorphic coefficients. (See Corollary 4.16 of \cite{OP04}.) Hence 
\begin{align*}
\frac{D\Psi_\e(p,X)}{p}&=\frac{D'(\usfrac{d}{dX})\Psi_\e(p,X)}{p}\\
&=\frac{2 \widehat{f}(ip-\rho) D'(\usfrac{d}{dX})\cos(pX)}{p}\\
           &=-2\widehat{f}(ip-\rho)D'\sin(pX)\,.
\end{align*} 
Formula (\ref{eq:duPsi}) then yields, for $b=\exp X \in B$ and $t \in 
[0,R-\e]$,
\renewcommand{\arraystretch}{2.3}
\begin{equation}\label{eq:durankone}
 \begin{array}{rl}
\delta(b)u(b,t)&\!\!\!=
D\left( \dint_{0}^\infty \dfrac{\Psi_\e(p,X)}{p}\, \sin(pt)\; 
dp \right) \\
&\!\!\!=
\dfrac{1}{2i} D\left( \dint_{-\infty}^\infty \dfrac{\Psi_\e(p,X)}{p}\,
  e^{ipt}\; 
dp \right) \\
&\!\!\!=
\dfrac{1}{2i} \dint_{-\infty}^\infty \dfrac{D\Psi_\e(p,X)}{p}\, e^{ipt}\; dp
 \\
&\!\!\!=
i \dint_{-\infty}^\infty \widehat{f}(ip-\rho) D'\sin(pX)\, e^{ipt}\; dp.
\end{array}
\end{equation}
\renewcommand{\arraystretch}{1.3}

Since $\widehat{f} \in {\rm PW}_\e(\fb^*)$, for all $N \in \N$ there are 
positive constants $C_N$ and $C_N'$ so that 
\begin{align*}
\big|\widehat{f}(ip-\rho) \sin(pX)\big|
&\leq C_N (1+\|ip-\rho\|)^{-N} e^{\e\|\re(ip-\rho)\|}
   e^{\re(ip\omega(X))}\\
&\leq C'_N (1+|p|)^{-N} e^{\e|\im p|} 
   e^{|\im p|\|X\|}\,
\end{align*}
for all $p\in \C$ and $X \in \fb$.
As in Lemma \ref{lemma:Psi}(b) we conclude that 
$\widehat{f}(ip-\rho) D'\sin(pX)$ is a holomorphic function of 
$(p,X) \in \C \times \fb$, and for every $N \in \N$ 
and every compact $Q \subset \fb$ there is a constant $K_{N,Q}>0$
such that
\begin{equation} \label{eq:estDPsirankone}
\Big|  \widehat{f}(ip-\rho) D'\sin(pX)\Big|
 \leq K_{N,Q}  (1+|p|)^{-N} 
e^{|\im p|(\e+\|X\|)}
\end{equation}
for all $p \in \C$ and $X \in Q$.
This allows us to shift the contour of integration in (\ref{eq:durankone})
and get for all $p\in \C$ and $b=\exp X \in B$:
 \begin{equation}\label{eq:shiftrankone}
\delta(b)u(b,t)=
\left(\dint_{-\infty}^\infty \widehat{f}(ip-\gamma+\rho) 
D'\sin((p+i\gamma)X)\, e^{ipt}\; dp \right) e^{-\gamma t}.
\end{equation}
The same argument used in Lemma \ref{lemma:expbigrank}, together with 
(\ref{eq:estDPsirankone}) and (\ref{eq:shiftrankone}), yields the
following lemma.

\begin{lemma} \label{lemma:exprankone}
Keep the assumptions of Theorem \ref{thm:HP}. Then
the local exponential Huygens' principle of Theorem \ref{thm:HP}(a)
holds when $\dim U/K =1$.
\end{lemma}

The local  exponential Huygens' principle provides a second proof of the
local strong Huygens' principle.

\begin{corollary} \label{lemma:strongHPdue}
Keep the assumptions of Theorem \ref{thm:HP}. If $\dim(U/K)$ is odd, then
the local strong Huygens' principle of Theorem \ref{thm:HP}(b)
holds for the modified wave equation
on $U/K$.
\end{corollary}
\begin{proof}
The finite propagation speed ensures that 
$$\Supp(u) \cap ( U/K \times [0,R-\e])= 
\Supp(u) \cap ( D_R \times [0,R-\e]) \subseteq C_\e\,,$$
where $C_\e$ denotes the positive $\e$-cone (\ref{eq:econe}).
As $\gamma \to \infty$, we obtain from (\ref{eq:expHP}) that 
$\delta(x)u(x,t)=0$ for all $x$ with 
$t-d(x,x_0)-\e>0$.
Thus $\Supp(u) \cap ( U/K \times [0,R-\e])$ is contained in the
$\e$-shell $S_\e$ of (\ref{eq:eshell}).
\end{proof}

We now turn to the proof of the explicit formulas for the smooth 
solution of the Cauchy problem (\ref{eq:cauchy}) for the 
modified wave equation on $U/K$. These formulas are
a consequence of (\ref{eq:deltau}) and of the explicit formulas 
known for the solution to the Cauchy problem (\ref{eq:euclCauchy})
for the Euclidean wave equation.

For $r>0$ we denote by $S_r(X):=\{X\in \fb\mid 
\|X\|=r\}$ the Euclidean sphere in $\fb \cong \R^n$ 
of center $X$ and radius $r$.
Again, we let $\Omega_{n-1}(r)$ denote the surface area of $S_r(x)$, and 
write simply $\Omega_{n-1}$ for $\Omega_{n-1}(1)$. Recall 
the definition (\ref{eq:M}) of the
mean value $(M^rf)(X)$ of 
a function $f:\fb\to \C$ on $S_r(X)$.

\begin{lemma}\label{lemma:solutionEucl}
Suppose $\dim \fb=n$ is at least $2$.

If $n$ is odd, then 
the solution to \straightparen{\ref{eq:euclCauchy}} is given by 
\begin{equation}
\label{eq:vodd}
v(X,t)=\frac{\Omega_n/2}{[(n-3)/2]! \Omega_{n-1}}\; 
\left( \dfrac{\partial}{\partial(t^2)}\right)^{(n-3)/2} 
  \big( t^{n-2} (M^tf)(X)\big)\,.
\end{equation}
If $n$ is even, then the solution to \straightparen{\ref{eq:euclCauchy}} 
is given by 
\begin{equation} \label{eq:veven}
v(X,t)=\frac{1/2}{[(n-2)/2]!} \int_0^t r (t^2-r^2) 
 \left(\dfrac{\partial}{\partial r^2}\right)^{(n-2)/2} \big(r^{n-2}
 (M^rf)(X)\big)
\; dr\,.
\end{equation}
\end{lemma}
\begin{proof}
See e.g.\ \cite{He94}, p.\ 481. 
\end{proof}
 
\begin{corollary}\label{cor:formulas}
Let $U/K$ be a symmetric space of the compact type with even multiplicities.
Suppose $\dim U/K \geq 2$. Let $D$ be the differential operator of
Theorem \ref{th-anaext}, and keep the notation of Lemma 
\ref{lemma:solutionEucl}.
If $n=\rank(U/K)=\dim \fb$ is odd, then 
the smooth solution $u(b,t)$ 
to \straightparen{\ref{eq:euclCauchy}} is given, for all $b=\exp X \in B$ 
and $t \in [0,R-\e]$, by the formula 
\begin{equation*}
\delta(b)u(b,t)=\frac{\Omega_n/2}{[(n-3)/2]! \Omega_{n-1}}\; 
D\left( \dfrac{\partial}{\partial(t^2)}\right)^{(n-3)/2} 
  \big( t^{n-2} (M^tf)(X)\big)\,.
\end{equation*}

If $n$ is even, then the solution to (\ref{eq:euclCauchy}) is given 
 for all $b=\exp X \in B$ 
and $t \in [0,R-\e]$ by the formula 
 \begin{equation}\label{eq:uoddbis} 
\delta(b)u(b,t)=\frac{1/2}{[(n-2)/2]!} D \int_0^t r (t^2-r^2) 
 \Big(\frac{\partial}{\partial r^2}\Big)^{(n-2)/2} \big(r^{n-2} (M^rf)(X)\big)
\; dr\,.
\end{equation}
\end{corollary}
\begin{proof}
This is immediate from (\ref{eq:deltau}), (\ref{eq:v}), and Lemma 
\ref{lemma:solutionEucl}.
\end{proof}

Corollary \ref{cor:formulas} proves, in particular, Theorem \ref{thm:HP}(c).
It also yields a third proof of the local strong Huygens' principle.
Indeed, (\ref{eq:uodd}) shows that $u(b,t)$ is determined by the values of the
Cauchy datum $f$ in a thin shell around $S_t(X)$, where $b=\exp X$.

\end{document}